\documentclass[11pt]{article}

\usepackage{amsthm,amsmath,amssymb,bm}
\usepackage{graphicx,color,epsfig}
\usepackage[margin=1in]{geometry}
\usepackage{verbatim}

\newtheorem{Thm}{Theorem}[section]
\newtheorem{Lem}[Thm]{Lemma}
\newtheorem{Pro}[Thm]{Proposition}
\newtheorem{Claim}[Thm]{Claim}

\newcommand{\D}{\mathbb{D}}
\newcommand{\E}{\mathbb{E}}
\newcommand{\N}{\mathbb{N}}

\newcommand{\R}{\mathbb{R}}

\newcommand{\Z}{\mathbb{Z}}

\begin{document}

\title{Slow-fast dynamics in stochastic Lotka-Volterra systems}
\author{Julien Barr\'e$^1$, Bastien Fernandez$^2$ and Gr\'egoire Panel$^1$}
\maketitle

\begin{center}

$^1$ Institut Denis Poisson, Universit\'e d'Orl\'eans, Universit\'e de Tours and CNRS, Orl\'eans, France.

$^2$ Laboratoire de Probabilit\'es, Statistique et Mod\'elisation\\
Univ.\ Paris Cit\'e -  CNRS - Sorbonne Univ.\\
F-75013 Paris, France
\end{center}

\begin{abstract}
We investigate the large population dynamics of a family of stochastic particle systems with three-state cyclic individual behaviour and parameter-dependent transition rates. On short time scales, the dynamics turns out to be approximated by an integrable Hamiltonian system whose phase space is foliated by periodic trajectories. 
%A fast dynamics takes place on short times scales, that is governed by an integrable Hamiltonian system whose phase space is foliated by periodic trajectories. 
This feature suggests to consider the effective dynamics of the long-term process that results from averaging over the rapid oscillations. We establish the convergence of this process in the large population limit to the solutions of an explicit stochastic differential equation. Remarkably, this averaging phenomenon is complemented by the convergence of stationary measures. The proof of averaging follows the Stroock-Varadhan approach to martingale problems and relies on a fine analysis of the system's dynamical features. 
\end{abstract}

\section{Introduction}
Stochastic interacting particle systems with cyclic structure, sometimes called stochastic Lotka-Volterra (LV) systems or rock-paper-scissors games, play an important role in modelling in a large variety of different fields: ecology and population dynamics \cite{BM15,BRSF09,CFM06,M10}, evolutionary game theory \cite{F10,SMJSRP14}, dynamics of Bose-Einstein condensates \cite{KWKF15}, chemical reaction networks  \cite{SNSWS17}, etc. To describe the long term dynamics of these models is therefore an important challenge in the theoretical and mathematical sciences. 

In the idealised large population limit, these stochastic systems are usually well-approximated by systems of ordinary differential equations (of LV type), 
see for instance \cite{AL84,F10,GKF18,KWKF15}; thus considerations about deterministic dynamics could suffice in principle. Yet, in the more realistic case of finite populations, the deterministic approximations are valid only on short time scales. Stochasticity must be integrated into the analysis in order to apprehend the long-term behaviors. In fact, dramatic finite-size effects can generate various phenomena, such as extinction events and other averaging phenomena,  which are not captured by the deterministic limit, see e.g.\ \cite{BRSF09,CT08,DF12,IP13,RMF06} for examples in the physics literature. From a rigorous mathematical viewpoint, extinction events have been described in examples of particle systems, in particular in population dynamics \cite{D08,E04}. However, to the best of our knowledge, no mathematical characterisation of an emerging averaging phenomena in stochastic LV systems has been given in the literature. Of note, averaging is ubiquitous in particle systems without cyclic behaviour, when the slow-fast time scale separation naturally materializes in the original variables. Various mathematical results have been obtained in this context, see e.g.\ multiscale chemical reaction and gene networks \cite{BKPR06,CDMR12,KK13,KKP14} and structured population dynamics \cite{C16,MT12}.

In order to mathematically address averaging in stochastic LV systems as it emerges from oscillatory behaviours that involve all the original variables, we consider in this paper a simple example of particle system with cyclic state space. In few words, the system is a Markov process that can be defined as follows (see section \ref{S-SYST} for details). The state space is $\Z_3=\Z/3\Z$ and a particle in state $i\in \Z_3$ can only jump to state $i+1$. The jumps are independent and for each particle in state $i$, occur with rate $a+N_{i+1}$ where $a\in\R^+$ is an intrinsic rate and $N_{i+1}\in\N$ is the number of particles in state $i+1$. The population size $N=\sum_{i\in\Z_3} N_i$ is constant so that the phase space is the two-dimensional simplicial grid. 

\begin{figure}[ht]
\begin{center}
\includegraphics*[width=60mm]{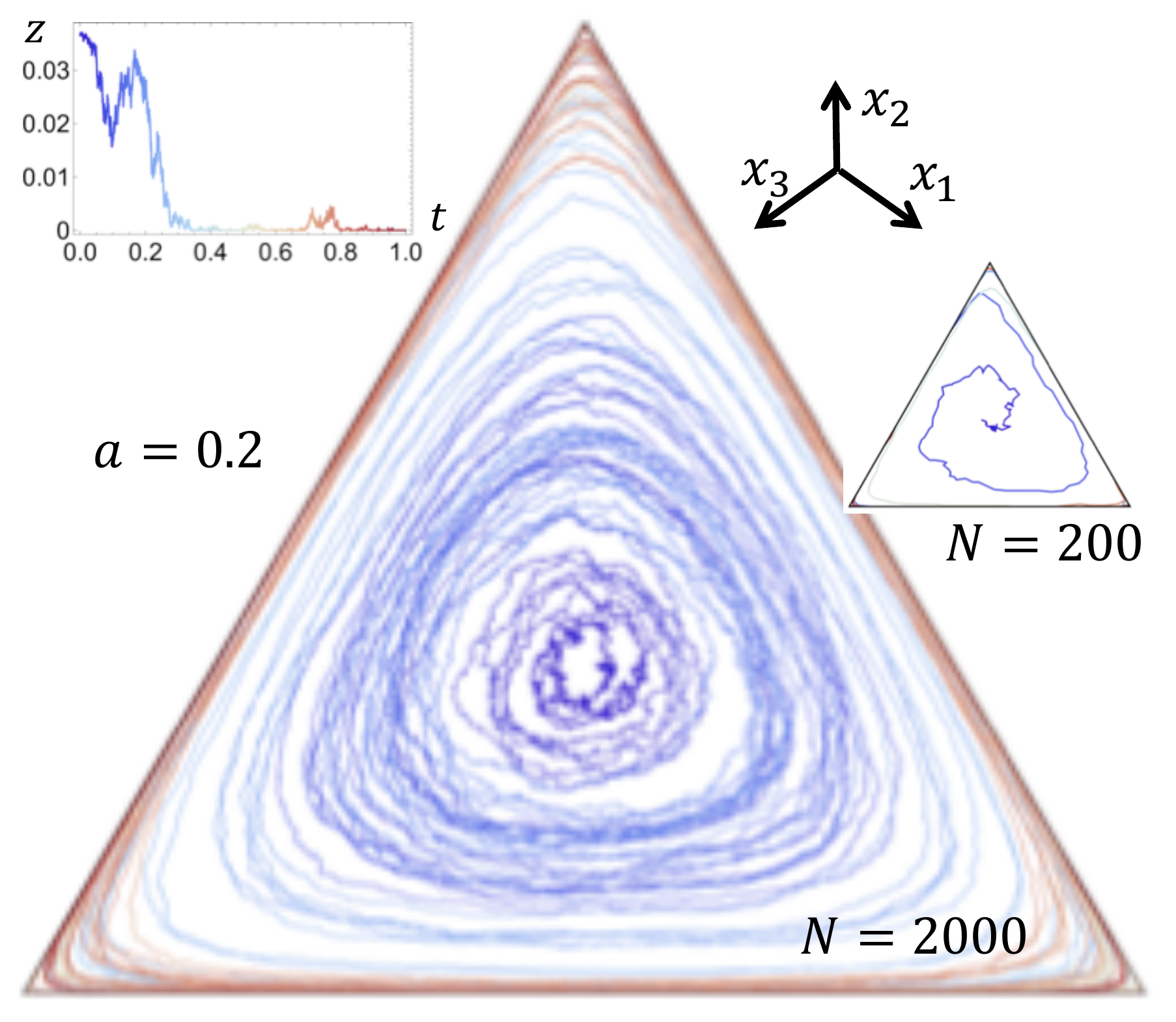}
\hspace{2cm}
\includegraphics*[width=60mm]{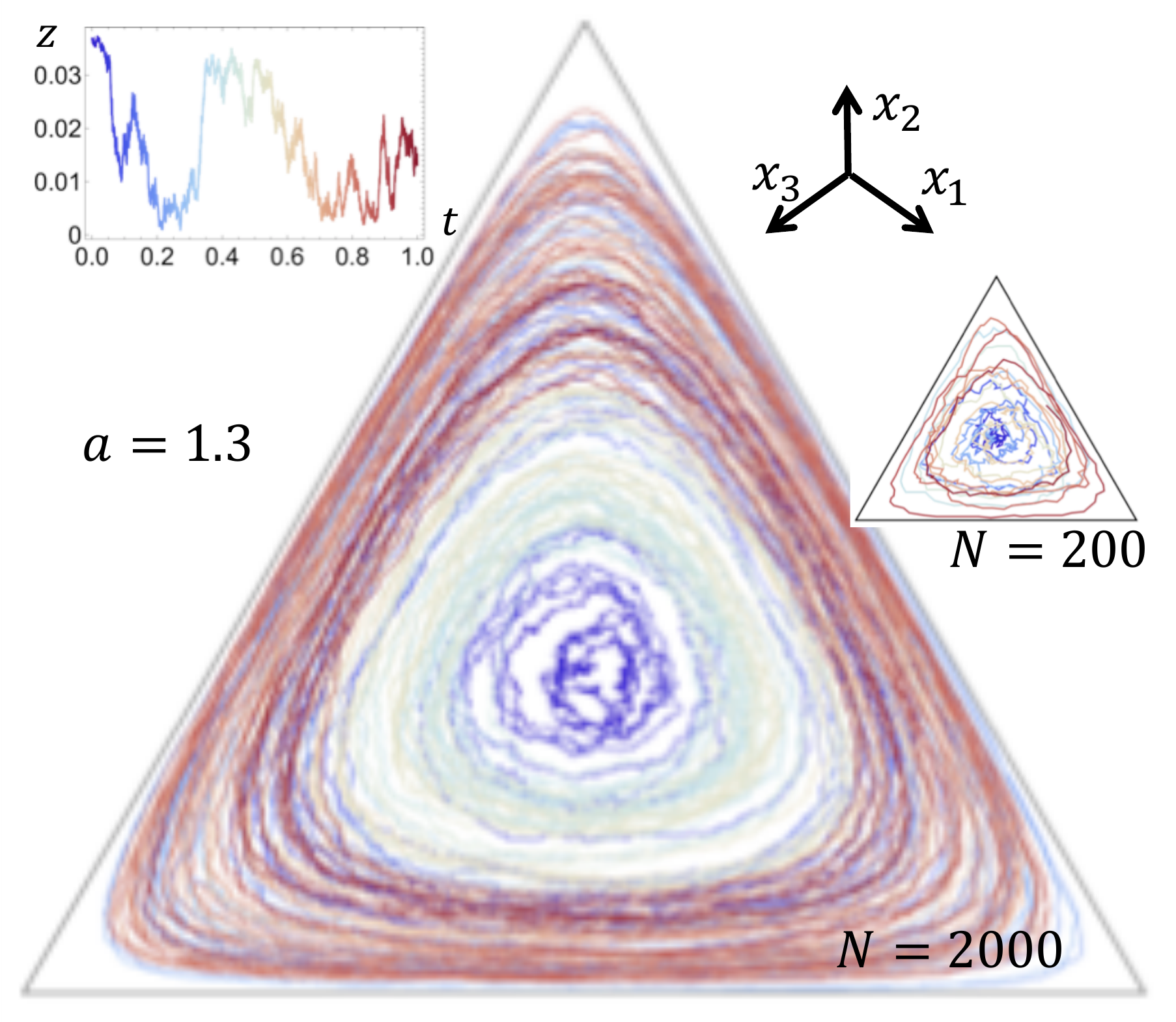}
\end{center}
\caption{Trajectories of the particle system in the two-dimensional simplicial grid ({\sl Left}\ $a=0.2$, {\sl Right}\ $a=1.3$, $N=2000$ in the main pictures, $N=200$ in the right insets). The initial condition are located at the center of the grid. The colors stand for the time in $[0,1]$, from blue ($t=0$) to red ($t=1$). {\sl Left inset:} Time series of the slow variable $z(t)=\prod_{i\in\Z_3}x_i(t)$ for $N=2000$, where $x_i=\tfrac{N_i(t)}{N}$.}
\label{NUMERICS}
\end{figure}
As intended, such extensive transition rates promote rapid oscillations in phase space when $N$ is large (see illustrations on Fig.\ \ref{NUMERICS}, in particular compare the main pictures $N=2000$ from the corresponding right insets $N=200$). In fact, the deterministic flow in the large population limit, which turns out to approximate the short time scale dynamics when $N$ is large (see Proposition \ref{APPROXFAST} below), consists of an integrable Hamiltonian system whose two-dimensional simplex phase space is foliated by periodic trajectories (on which the slow variable $z=\prod_{i\in\Z_3}x_i$ remains constant). This suggests to consider the one-dimensional transverse dynamics of the variable $z$ that results from averaging the fast motions on the periodic loops. The main result of this paper (Theorem \ref{MAINRES}) states that for large $N$, the slow-time scale transverse dynamics is indeed approximated by a diffusion process with $a$-dependent drift. In short terms, a slow-fast dynamics emerges in this system in the large population limit. 

Technically speaking, the stochastic process that governs the dynamics of the particle system. can be regarded as a random perturbation of a dynamical system with a conservation law. Yet, the oscillation period diverges at the phase space boundary (independently of the population size) and this prevents us to apply the standard techniques in this setting \cite{FW04,PS08}. Instead, our proof follows the Stroock-Varadhan approach to martingale problems \cite{SV79} and relies on the compactness-uniqueness argument in this context. The core argument (section \ref{S-INDENT}) is a proof of the $L^1$-convergence of martingales which is tailored to the specific nature of the process and in particular, to its behaviour close to the boundary. 

Remarkably, for this particle system, the averaging phenomenon is further complemented by the large $N$ convergence of stationary measures. Indeed, for every $N$, the (unique) stationary measure of the process on the simplicial grid is a product measure which converges to a Dirichlet distribution in the large population limit (Proposition \ref{MU_INV_N_INFINI}). Moreover, the push-forward measure on the transverse variable $z$ induced by this distribution turns out to be stationary for the semi-group associated with the diffusion process (Proposition \ref{INVARIANT_MEASURE}). Together with the specification of the nature of the boundary points of this process (Lemma \ref{pro:boundary}), these properties indicate that the visits of the particles' system to the boundaries of the simplex are frequent for $a<1$ and become sparse when $a\geq 1$, as illustrated in Fig.\ \ref{NUMERICS}. 

Our results are limited here to a simple model with three states playing symmetric roles, but the ideas and techniques can be useful in a broader context. For instance, the analysis can be extended to state-dependent transition rates, which may provide clues to answer the important question: when extinction is possible, which species survives? The answer is counter intuitive as shown in the physics literature  \cite{BRSF09}. The extension to more than three states is more challenging, as the deterministic dynamics may not be periodic any more; however, when it is periodic, the techniques of this article may allow to prove the convergence of the slow dynamics to a multi dimensional diffusion process.

\section{Definitions and preliminary considerations}
\subsection{The stochastic particle system}\label{S-SYST}
\begin{figure}[ht]
\begin{center}
\includegraphics*[width=60mm]{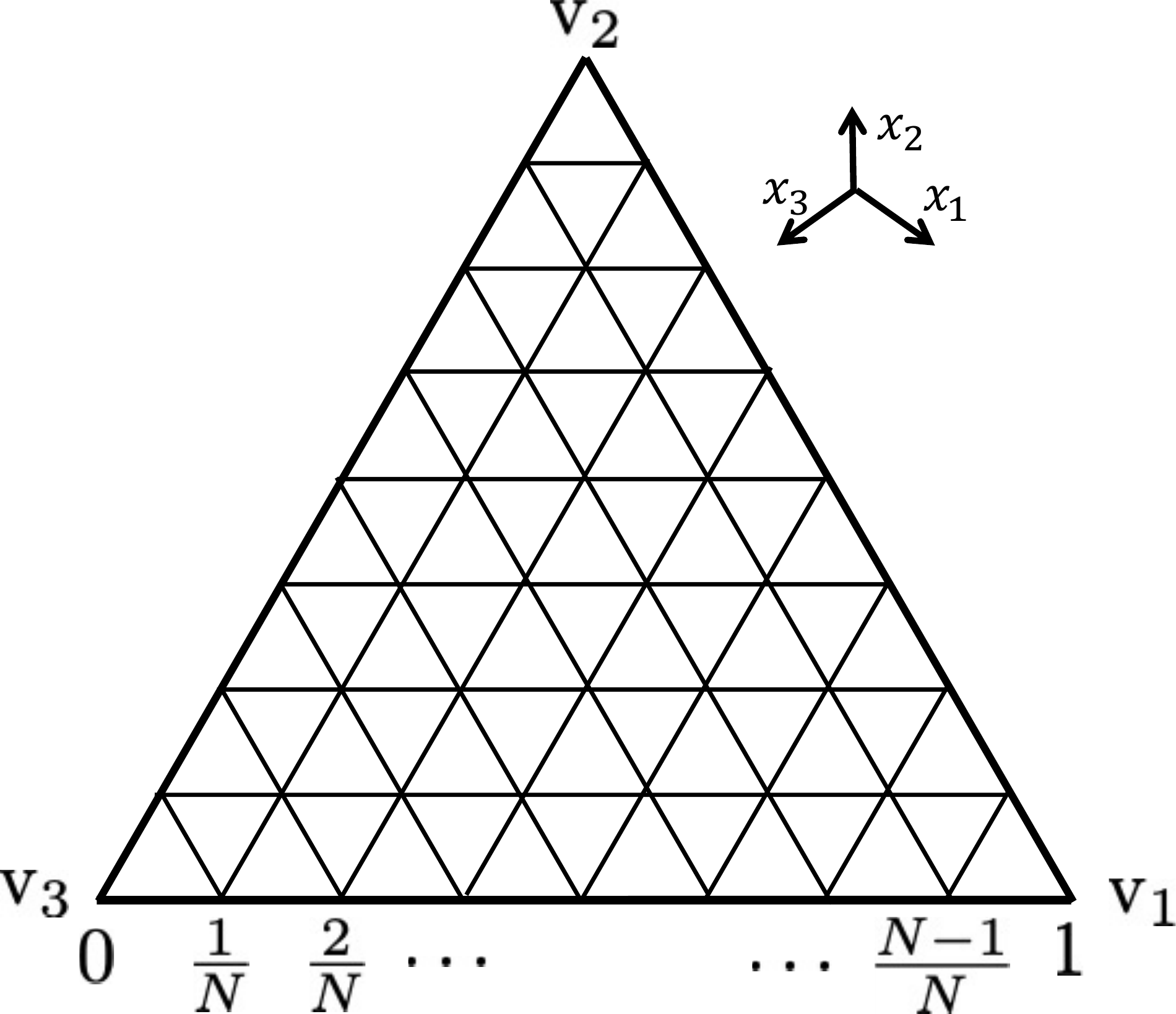}
\end{center}
\caption{Illustration of the simplicial grid $S_N$ of step size $\tfrac1{N}$.}
\label{SIMPLGRID}
\end{figure}

We consider the two-dimensional simplex $S$ defined by 
\[
S=\left\{\mathrm{x}:=(x_1,x_2)\in (\R^+)^2\ \text{such that}\ x_1+x_2\leq 1\right\},
\]
and given $N\in\N$ (where $\N=\{0,1,2,\cdots\}$), let $S_N=S\cap \tfrac1{N}\N^2$ be the two-dimensional simplicial grid, whose vertices $\mathrm{v}_i$ are the points with coordinates $(\mathrm{v}_i)_j=\delta_{ij},\,j=1,2$ (where $\delta_{ij}$ is the Kronecker symbol)  (see Fig.\ \ref{SIMPLGRID}).

The time evolution of the particle system in $S_N$ is governed by the (jump) Markov process $\{P_N^\mathrm{x}\}_{\mathrm{x}\in S_N}$ induced by the generator $L_N$ defined by\footnote{Throughout the paper, the notations $i+1$ and $i-1$ mean respectively $i+1\ \text{mod}\ 3$ and $i-1\ \text{mod}\ 3$.}  
\[
L_Nf(\mathrm{x})=N\sum_{i\in\Z_3}x_i(a+Nx_{i+1})\left(f(\mathrm{x}+\frac{\mathrm{u}_i}{N})-f(\mathrm{x})\right)\quad \forall \mathrm{x}\in S_N,\ f:S_N\to \R
\]
where $x_3:=1-x_1-x_2$, $a\in\R^+$ and $\mathrm{u}_i=v_{i+1}-v_i$ for $i\in\Z_3:=\Z/3\Z$.

As mentioned in the introduction, this process represents the stochastic time evolution of a population of individuals with cyclic state space and extensive transition rates and is inspired by the modelling in various fields \cite{FT08,RMF06}. In particular, the definition above suggests various natural extensions of this process, such as increasing the number of states, from three to an arbitrary $d\in\N$, or allowing any particle of a site $i$ to jump on the site $j$. Notice that most of the approach and considerations in this paper can be adapted to these extensions without additional conceptual difficulties.

A nice feature of this process is that it is ergodic for every $N\in\N$ and $a>0$ and its invariant measure turns out to be the following product measure \cite{FT08}
\[
\mu_{N,a}(\mathrm{x}) = C_{N,a} \prod_{i\in\Z_3}\frac{\Gamma(Nx_i+a)}{\Gamma(Nx_i+1)},\quad \forall \mathrm{x}\in S_N,
\] 
where $\Gamma$ stands for the Gamma function and $C_{N,a}$ is the normalisation constant. For $a=0$, the three vertices $\{\mathrm{v}_i\}_{i\in\Z_3}$ are absorbing states. 

The measure $\mu_{N,a}$ can be regarded as an atomic measure in $S$. Under this viewpoint, this measure can be shown to weakly converge to the Dirichlet measure, namely the absolutely continuous measure $\mu_a$ on $S$ with density $\rho_a$ defined by 
\[
\rho_a(\mathrm{x})=C_a z(\mathrm{x})^{a-1},\ \forall \mathrm{x}\in S\ \text{where}\ z(\mathrm{x})=\prod_{i\in\Z_3}x_i,
\]
and again $C_a$ is the normalisation constant. The convergence is claimed in the following statement, whose proof is given in Appendix \ref{A-MU_INV_N_INFINI}.
\begin{Pro}
For every $a>0$, we have
\[
\lim_{N\to\infty} \mu_{N,a}=\mu_a,
\]
in the weak sense. 
\label{MU_INV_N_INFINI}
\end{Pro}

\subsection{The deterministic approximation on short time scales}
When interested in the temporal process associated with $L_N$ for large $N$, the Taylor theorem applied to the first order expansion of $f(\mathrm{x}+\frac{\mathrm{u}_i}{N})$ at $\mathrm{x}$ suggests to consider the operator ${\cal L}_\text{fast}$ defined by (NB: $f'$ denotes the Fr\'echet derivative of $f$.)
\[
{\cal L}_\text{fast}f(\mathrm{x})=\sum_{i\in\Z_3}x_ix_{i+1}f'(\mathrm{x})\mathrm{u}_i =\sum_{j=1,2}x_j(x_{j-1}-x_{j+1})\partial_{x_j}f(\mathrm{x}),\ f\in C^1(S),
\]
so that we have $\lim_{N\to\infty}\tfrac1{N}L_Nf(\mathrm{x})= {\cal L}_\text{fast}f(\mathrm{x})$. 
Let $F$ be the vector field on $S$ defined by 
\[
(F(\mathrm{x}))_j=x_j(x_{j-1}-x_{j+1}),\ j=1,2.
\]
This vector field defines a semi-flow on $S$, under which this simplex is invariant.
Let $t\mapsto X^{\rm fast}_{\mathrm{x}_0}(t)$ be the solution of $\dot{\mathrm{x}}=F(\mathrm{x})$ with initial condition $\mathrm{x}_0\in S$. Then for any $f\in C^1(S)$, we have
\[
\tfrac{d}{dt}f(X^{\rm fast}_{\mathrm{x}_0}(t))={\cal L}_\text{fast} f(X^{\rm fast}_{\mathrm{x}_0}(t)),\ t>0.
\]
The convergence $\tfrac1{N}L_Nf(\mathrm{x})\to {\cal L}_\text{fast}f(\mathrm{x})$ suggests that the process associated with the particle system can be approximated on the time scales of the order $\frac1{N}$ by the deterministic semi-flow. In order to formalize this approximation, given $T>0$, let $\D([0,T],S_N)$ (resp.\ $\D([0,T],S)$) be the set of c\`adl\`ag functions from $[0,T]$ into $S_N$ (resp.\ $S$) and let ${\cal F}_{T,N}$ (resp.\ ${\cal F}_{T}$) be the natural filtration associated with $\D([0,T],S_N)$ (resp.\ $\D([0,T],S)$). The set $\D([0,T],S)$ is endowed with the Skorokhod metric. We denote by $X_N(t)$ where $t\in [0,T]$ and $X_N\in \D([0,T],S_N)$ the stochastic process on $(\D([0,T],S_N),{\cal F}_{T,N})$ associated with the Markov process and the initial measure $\delta_{X_N(0)}$. Clearly, $X_N(t)$ can be seen as a process taking values in $S$ 
\begin{Pro}
Assume that the sequence of initial conditions $\{X_N(0)\}_{N\in\N}$ converges in law to some $\mathrm{x}_0\in S$. Then, for every $T>0$, the sequence of time-scaled processes $\{X_N(\frac{t}{N}) : t\in [0,T]\}_{N\in\N}$ converges in law in $\D([0,T],S)$ to the trajectory arc $\{X^{\rm fast}_{\mathrm{x}_0}(t) : t\in [0,T]\}$ of the solution of $\dot{\mathrm{x}}=F(\mathrm{x})$ with initial condition $\mathrm{x}_0$.
\label{APPROXFAST}
\end{Pro} 
\noindent
This statement can be proved using a compactness-uniqueness argument just as in the proof of Theorem 3.1, Chap.\ 3 in \cite{BM15}. See also the proof of Theorem \ref{MAINRES} below for the details of a similar argument. 

\subsection{Analysis of the deterministic dynamics}
According to the expression of $F$, the semi-flow associated with $\dot{\mathrm{x}}=F(\mathrm{x})$ is an instance of a Lotka-Volterra system. 
Actually, this dynamics can be regarded as a Hamiltonian system with Hamiltonian function $(x_1,x_2)\mapsto x_1x_2(1-x_1-x_2)$.
%\cite{I87,M75,TKW97}. 
%Moreover, since the functional $z(\mathrm{x})$ is invariant, this system is also an elementary instance of an integrable Hamiltonian system.
% \cite{B08,DEKV17}. 

The dynamics can be analysed in full details and its essential features have already been identified \cite{AL84,RMF06}. In particular, there are four stationary points, namely the centre $(\tfrac13,\tfrac13,\tfrac13)$ and the vertices $\mathrm{v}_i$ of $S$. Each boundary edge of $S$ is invariant under the semi-flow and the dynamics on each edge consists of heteroclinic trajectories between the two corresponding vertices. 

In addition, in $\text{Int}(S)\setminus (\tfrac13,\tfrac13,\tfrac13)$, the interior of the simplex except the centre, the level sets of the functional $z$ - which takes values in $I:=(0,\tfrac1{27})$ - constitute a foliation by invariant loops on which the trajectories $t\mapsto X^{\rm fast}_{\mathrm{x}_0}(t)$ are periodic, with period say $T(z(\mathrm{x}_0))$, and counterclockwise motion (see Fig.\ \ref{CONTOURS}). The periodic trajectories and their period can be semi-explicitly computed, see Appendix \ref{A-DYNAML0-1} for the corresponding computations. In particular, the period diverges when approaching the boundary edges of $S$. 
\begin{figure}[ht]
\begin{center}
\includegraphics*[width=70mm]{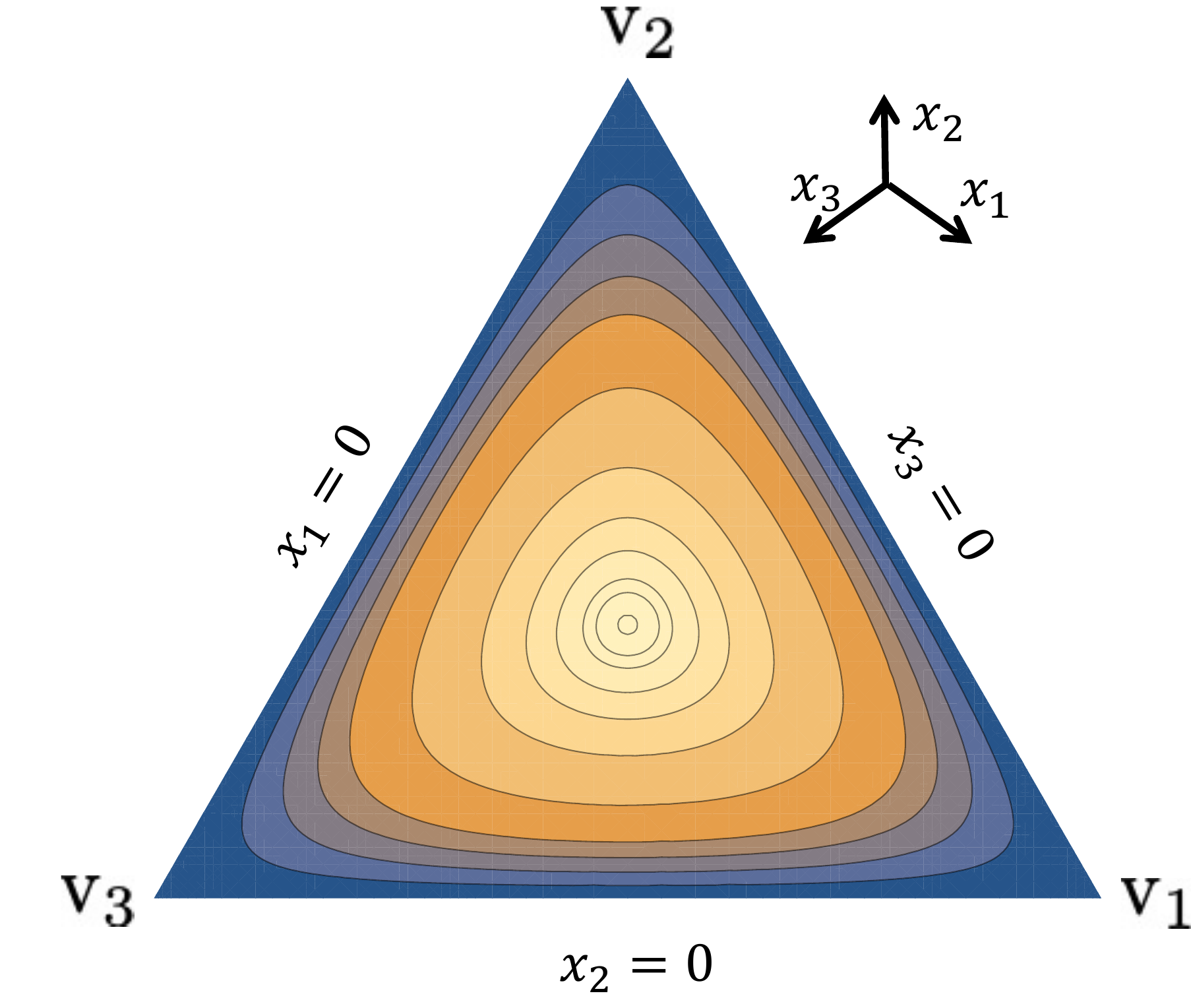}
\end{center}
\caption{Color plot and level sets of the function $\mathrm{x}\mapsto z(\mathrm{x})$ in the simplex $S$.}
\label{CONTOURS}
\end{figure}

In the sequel, we shall need the following additional properties of the period function. Of note, we use the symbol $z$ for the variable in $I$ and also $T(z)$ as an abbreviation of the notation of the period. Moreover, given two real functions $u$ and $v\neq 0$ and $x_0\in\R$, we write $u(x)\sim v(x)$ as $x\to x_0^\pm$ if $\lim_{x\to x_0^\pm}\frac{u(x)}{v(x)}=1$.
\begin{Lem}
(i) The function $z\mapsto T(z)$ is $C^\infty$ on $I$ and $T(\frac1{27}^-)=2\pi\sqrt{3}$.

\noindent
(ii) We have $T(z) \sim -3\ln z$ as $z\to 0^+$. 
\label{PROPERIOD}
\end{Lem}
\noindent
The proof is given in Appendix \ref{A-PROPERIOD}. In addition, we shall also need some properties of the (signed) area enclosed in the loop $z(\mathrm{x})=z$ and defined by\footnote{Lemma \ref{AVGOP} below shows that in fact $A(z)=\int_0^{T(z)}x_{i+1}\dot{x}_idt$ for every $i\in\Z_3$.}
\[
A(z)=\int_0^{T(z)}x_2\dot{x}_1dt,\quad z\in I,
\]
which can be regarded as the action of the Hamiltonian system. The desired properties of this function are listed in the next statement, whose proof is given in Appendix \ref{DYNRAP-A}.
\begin{Lem}
(i) The function $z\mapsto A(z)$ is $C^\infty$ and negative on $I$. Moreover, we have $A'(z)=T(z)$ for all $z\in I$. 

\noindent
(ii) $A(0^+)=-\frac12$ and $A(z)\sim -2\pi\sqrt{3}(\frac1{27}-z)$ as $z\to \frac1{27}^-$. 
\label{pro:dyn_rapide}
\end{Lem}

\subsection{The averaged generator: definition and explicit expression} 
Following the approach to Anosov averaging \cite{PS08}, the Proposition \ref{APPROXFAST} and the periodic motions of the system $\dot{\mathrm{x}}=F(\mathrm{x})$ suggest to consider averaging the dynamics associated with the next-order approximation of $L_N$. To that goal, consider first the operator ${\cal L}_\text{slow}$ that collects the $N$-independent terms in the (second-order) expansion of $L_N$, and defined for $f\in C^2(S)$ as follows
\begin{align*}
{\cal L}_\text{slow}f(\mathrm{x})&=a\sum_{i\in\Z_3}x_if'(\mathrm{x})\mathrm{u}_i+\tfrac12\sum_{i\in\Z_3}x_ix_{i+1}f''(\mathrm{x})(\mathrm{u}_i,\mathrm{u}_i)\\
&=a\sum_{j=1,2}(x_{j+1}-x_j)\partial_{x_j} f(\mathrm{x})+ \tfrac12\sum_{j=1,2}x_j(x_{j-1}+x_{j+1})\partial_{x_j}^2f(\mathrm{x})-\sum_{j=1,2}x_jx_{j+1}\partial_{x_j,x_{j+1}}^2f(\mathrm{x})
\end{align*}
Given $z\in I$ and a fonction $f$ defined on the loop of period $T(z)$, let the time average $\langle f\rangle_z$ be defined by
\[
\langle f\rangle_z=\frac1{T(z)}\int_0^{T(z)}f(X^{\rm fast}_{\mathrm{x}_0}(t))dt,
\]
(where $\mathrm{x}_0$ is any point on the loop). Then, using the notation $f_I$ for the functions defined on $I$ (or $\bar{I}$), the averaged operator ${\cal L}_\text{avg}$ is defined by
\[
{\cal L}_\text{avg}f_I(z)=\langle {\cal L}_\text{slow}(f_I\circ z)\rangle_z,\ z\in I,\ f_I\in C^2(I).
\]
An explicit expression of this operator can be obtained based on the analysis of the dynamics generated by $F$. The results are summarized in the following statement.
\begin{Lem}
The average $\langle x_{i+1}\dot x_i\rangle_z$ does not depend on $i\in\Z_3$. Letting $m(z):=-\langle x_{i+1}\dot x_i\rangle_z$, the averaged operator can be expressed as the following second-order differential operator
\[
{\cal L}_\text{\rm avg}f_I(z)=3(a m(z)-z)f_I'(z)+3z m(z)f_I''(z),\quad z\in I,\ f_I\in C^2(I).
\]
\label{AVGOP}
\end{Lem}
\noindent
Given the definition of $A(z)$ above, we have $m(z)=-\frac{A(z)}{T(z)}> 0$ for all $z\in I$ and Lemma \ref{pro:dyn_rapide} implies that $m \in C^\infty(I)$, which yields that  ${\cal L}_\text{avg}f_I \in C^0(I)$. Moreover, Lemma \ref{pro:dyn_rapide} and Lemma \ref{AVGOP} both imply that for $f_I\in C^2(\overline{I})$; hence ${\cal L}_\text{avg}f$ can be extended by continuity to the boundary of $I$ as follows
\[
{\cal L}_\text{avg}f_I(0) = 0\quad\text{and}\quad {\cal L}_\text{avg}f_I(\tfrac{1}{27})=-\tfrac19f_I'(\tfrac1{27}),
\]
which in particular yields the following match
\begin{equation}
{\cal L}_\text{avg}f_I(\tfrac{1}{27})={\cal L}_\text{slow}(f_I\circ z)(\tfrac13,\tfrac13,\tfrac13).
\label{MATCH}
\end{equation}
\noindent
{\sl Proof of the Lemma.} For $f_I\in C^1(I)$ and $\mathrm{x}\in \text{Int}(S)$, direct computations yield the following expression
\begin{align*}
{\cal L}_\text{slow}(f_I\circ z)&=\Big(a\sum_{i\in\Z_3}\big(x^2_{i-1}x_{i+1}-z\big)-3z\Big)f_I'\circ z+\tfrac{1}{2}\sum_{i\in\Z_3}z\Big(x_ix^2_{i-1}+x_ix^2_{i+1}-2z\Big)f_I''\circ z
\end{align*}
Together with the relations $z(\mathrm{x})\tfrac{x_{i-1}}{x_i}=x_{i-1}^2x_{i+1}$ and $z(\mathrm{x})\tfrac{x_{i+1}}{x_i}=x_{i-1}x_{i+1}^2$, this suggests consider the averages $\langle x_ix_{i+1}^2\rangle_z$ and $\langle x_ix_{i-1}^2\rangle_z$ in order to compute of the expression of ${\cal L}_\text{avg}$. We have the following statement.
\begin{Claim}
For every $z\in I$, we have $\langle x_ix_{i+1}^2\rangle_z=\langle x_ix_{i-1}^2\rangle_z$ for all $i\in\Z_3$ and these quantities do not depend on $i$.
\end{Claim}
\noindent
{\sl Proof:} By periodicity of the trajectories, we have for every $f\in C^1(S)$ and $z\in I$
\[
0=\langle \tfrac{d}{dt}f\rangle_z=\langle {\cal L}_\text{fast}f\rangle_z=\langle \sum_{j=1,2}x_j(x_{j-1}-x_{j+1})\partial_{x_j}f\rangle_z
\]
In particular, for $f(\mathrm{x})=x_ix_{i+1}$ for some $i\in \Z_3$, we get $\langle x_{i-1}x_{i+1}^2-x_{i+1}x_{i-1}^2\rangle_z=0$, which immediately yields the desired equality. In order to prove that the quantities do not depend on $i$, apply the equality above with $f(\mathrm{x})=x_{i+1}$, which combined with the previous one yields $\langle x_ix_{i+1}^2\rangle_z=\langle x_{i+1}x_{i-1}^2\rangle_z$. \hfill $\Box$

The expression of ${\cal L}_\text{avg}$ then immediately follows from the relation 
\[
\langle x_ix_{i+1}^2\rangle_z-z=\langle x_ix_{i+1}(x_{i+1}-x_{i-1})\rangle_z=-\langle x_{i+1}\dot x_i\rangle_z
\]
Lemma \ref{AVGOP} is proved. \hfill $\Box$

\subsection{The stochastic differential equation and its solutions}
Consider the (one-dimensional) stochastic differential equation associated with ${\cal L}_\text{avg}$, namely 
\begin{equation}
dZ(t)=b(Z(t))dt+\sigma(Z(t))dW(t),\quad Z(0)=z_0\in I \label{eq:eds_avg}
\end{equation}
where $b(z)=3(a m(z)-z)$, $\sigma(z)=\sqrt{6z m(z)}$ and $W(t)$ is some Brownian motion. Lemma \ref{pro:dyn_rapide} and the fact that $m(z)=-\frac{A(z)}{T(z)}$ for all $z\in I$ imply that both functions $b$ and $\sigma$ are smooth on $I$ and $\sigma^2>0$. These conditions ensure the existence and uniqueness of a solution $Z_{z_0}(t)$ for $t$ up to the so-called explosion time $T_\text{ex}(z_0)$, namely the time it takes for the solution to reach the boundary of $I$, see for instance \cite{IW89,KS91}. Of course, we have $T_\text{ex}>0$ a.s.

The explosion time depends on the nature of the boundary points, which can be evaluated using Feller's test \cite{IW89,KS91}. This nature depends on the parameter $a$ and is given in the following statement, which uses the Feller's classification in chapter 8.1 of \cite{EK86}. 
\begin{Lem}
\label{pro:boundary}
In the SDE \eqref{eq:eds_avg}, the boundary point $z=\frac{1}{27}$ is entrance for all $a\in \R^+$. Moreover the boundary point $z=0$ is entrance for $a\geq 1$, regular for $a\in (0,1)$ and exit for $a=0$.
\end{Lem}
\noindent
{\sl Proof:} We follow the arguments in chapter 8.1 of \cite{EK86}. Using $A'(z)=T(z)$ in Lemma \ref{pro:dyn_rapide}, the averaged generator can be recast as
\begin{equation}
\mathcal{L}_{\rm avg} =\tfrac{d}{ds(z)}\left(\tfrac{d}{d p(z)}\right),
\label{Lavg_scale_speed}
\end{equation}
where the scale function $p$ and the speed function $s$ are the positive functions whose differential are respectively given by the following equalities
\[
dp(z)=-\frac{1}{z^a A(z)}dz\quad \text{and}\quad ds(z)=\frac{z^{a-1}T(z)}3dz,\quad z\in I.
\]
The approximations in Lemmas \ref{PROPERIOD} and \ref{pro:dyn_rapide} imply the following ones
\[
\tfrac{dp(z)}{dz}\sim 2 z^{-a}\quad \text{and}\quad \tfrac{ds(z)}{dz}\sim - z^{a-1}\ln z\quad \text{as}\ z\to 0^+,
\]
and 
\[
\tfrac{dp(z)}{dz}\sim \tfrac{27^a}{2\pi\sqrt{3}}\left(\tfrac1{27}-z\right)^{-1}\quad \text{and}\quad \tfrac{ds(z)}{dz}\sim \tfrac{18\pi\sqrt{3}}{27^a}\quad \text{as}\ z\to \tfrac1{27}^-.
\]
Explicit computations then yield the following estimates for every $r\in I$ (recall that $a\in\R^+$)
\[
\left|\int_r^{0} s(z) dp(z)\right| <+\infty\ \text{iff}\ a\in [0,1)\quad \text{and}\quad \left|\int_r^0 p(z)ds(z)\right| <+\infty\ \text{iff}\ a>0,
\]
and
\[
\int_r^{\frac{1}{27}} s(z) dp(z) =+\infty\quad \text{and}\quad \int_r^{\frac{1}{27}} p(z)ds(z) <+\infty,
\]
from where the characterization of the boundary points $z=0$ and $z=\frac1{27}$ immediately follow. \hfill $\Box$
\medskip

As a consequence of the Lemma, for $a\geq 1$, we have $T_\text{\rm ex}=+\infty$ a.s., and hence existence and uniqueness of solutions of the SDE \eqref{eq:eds_avg} for all $t>0$, a.s. For $a\in [0,1)$, the trajectory hits the boundary point $z=0$ a.s. Yet, for $a=0$, existence and uniqueness of solutions for all $t>0$ a.s.\ follows by letting $Z(t)=0$ for $t>T_{\rm ex}$. 

For $a\in (0,1)$, the existence of solutions of the SDE \eqref{eq:eds_avg} extends for $t$ beyond $T_\text{ex}$ but uniqueness is not granted in general and requires to specify the behavior at the boundary point $z=0$ \cite{H96}. These features, as well as those for $a>1$ can be expressed through the semi-group associated with ${\cal L}_\text{avg}$. Following the definitions in \cite{EK86}, given $a>0$ consider the domain $\mathcal{D}_\text{avg}$ defined by 
\[
\mathcal{D}_\text{avg}=\left\{\begin{array}{ccl}
\left\{f_I \in C^2(I)\cap C^0(\bar{I})~:~ \mathcal{L}_{\rm avg}f_I \in  C^0(\bar{I}) \right\}&\text{if}&a\geq 1,\\
\left\{f_I \in C^2(I)\cap C^0(\bar{I})~:~ \mathcal{L}_{\rm avg}f_I \in  C^0(\bar{I}),~\lim_{z\to 0^+} z^aA(z)f_I'(z)=0 \right\}&\text{if}&a\in (0,1).
\end{array}\right.
\]
In particular, the choice of $\mathcal{D}_\text{avg}$ for $a\in (0,1)$ corresponds to an instantaneous reflexion at $z=0$. Theorem 1.1, Chap.\ 8 in \cite{EK86} states that, with these definitions of $\mathcal{D}_\text{avg}$, the operator ${\cal L}_\text{avg}$ generates a Feller semi-group on $C^0(\overline{I})$, for every $a>0$ (NB: In particular, in the case $a\in (0,1)$ for which $z=0$ is regular, the set $\mathcal{D}_0$ in \cite{EK86} is defined with $q_0=0$.)

Our next result states that the measure on $\overline{I}$ with density $z\mapsto z^{a-1} T(z)$ - which, thanks to the equality $A'(z)=T(z)$, results to be the push-forward under $\mathrm{x}\mapsto z(\mathrm{x})$ of the limit measure $\mu_a$ in Proposition \ref{MU_INV_N_INFINI} - is a stationary measure of the Feller semi-group (NB: Recall that for $a=0$, the point $z=0$ is exit). 
\begin{Pro} \label{INVARIANT_MEASURE} For every $a>0$, we have
\[
\int_0^{\frac{1}{27}} z^{a-1} T(z) \mathcal{L}_{\rm avg}f_I (z) dz =0,\ \forall  f_I \in \mathcal{D}_{\rm avg}.
\]
\end{Pro}
\noindent
{\sl Proof:} For $z_0<z_1\in I$, we obtain after direct integration
\[
\int_{z_0}^{z_1} z^{a-1} T(z) \mathcal{L}_{\rm avg}f_I (z) dz=\left[\tfrac{df_I(z)}{dp(z)}\right]_{z_0}^{z_1}.
\]
Therefore, all we have to show is $\lim_{z\to 0^+}\tfrac{df_I(z)}{dp(z)}=\lim_{z\to \frac1{27}^-}\tfrac{df_I(z)}{dp(z)}=0$. We focus on the first limit; the second one follows from similar argument. When $z=0$ is a regular boundary ($a\in (0,1)$), this is a consequence of the choice of $\mathcal{D}_{\text{avg}}$. When $z=0$ is entrance ($a\geq 1$), the Taylor formula implies that we have 
\[
\tfrac{df_I(z)}{dp(z)}=\frac{1}{p(z)-p(r)}\left(f_I(z)-f_I(r) +\int_{r}^{z}\,\mathcal{L}_{\rm avg}f_I(z)\,p(z)\,\mathrm{d}s(z)\right),\ r\in (z,\tfrac1{27}).
\]
The term inside the brackets converges when $z\to 0^+$ because we have $f_I,\mathcal{L}_{\rm avg}f_I \in C^0(\bar{I})$ and $|\int_{r}^{0}p(z) ds(z)|<+\infty$.  The result then follows from the limit $\lim_{z\to 0^+}|p(z)|=+\infty$. \hfill $\Box$

\section{Main result: Averaging}
Given $T>0$, let $\D([0,T],\overline{I})$ be the set of c\`adl\`ag functions from $[0,T]$ into $\overline{I}=[0,\frac1{27}]$ and let ${\cal G}_T$ be the natural filtration associated with $\D([0,T],\overline{I})$. The process $X_N(t)$ induces a stochastic process $Z_N(t):=z(X_N(t))$ on $(\D([0,T],\overline{I}),{\cal G}_T)$.
We are now in position to formulate the main result of the paper. 
\begin{Thm}
Assume that the sequence of initial conditions $\{Z_N(0)\}_{N\in\N}$ converges in law to some $z_0\in I$. \\
If $a\geq 1$ or $a=0$, then for any $T>0$, the sequence of processes $\{Z_N(t) : t\in [0,T]\}_{N\in\N}$ converges in law to the 
weak solution of the SDE \eqref{eq:eds_avg} with initial condition $z_0$. \\
If $0<a<1$, then for any $T>0$, the sequence $\{Z_N(t) : t\in [0,T]\}_{N\in\N}$ is relatively compact. Furthermore, any limit point is a weak solution of the SDE \eqref{eq:eds_avg}, with initial condition $z_0$.
\label{MAINRES}
\end{Thm}
\noindent
The proof of this statement follows the Stroock-Varadhan approach to martingale problems \cite{SV79}. The first step is a standard compactness argument (see the presentation in \cite{JM86}, especially Corollary 2.3.3 therein) for the semi-martingale structure associated with a stochastic process, the process $Z_N$ in our case. Then we prove that every limit point of a subsequence must satisfy the martingale problem associated with the SDE. That part of the proof is specific to the particle system under consideration as it relies in particular on various features of the short times deterministic dynamics. The convergence follows suit when the SDE has a unique solution, ie.\ for $a\geq 1$ and $a=0$. In the other case ($a\in (0,1)$), while the SDE solutions are not unique beyond $T_\text{ex}$, Proposition \ref{MU_INV_N_INFINI} and Proposition \ref{INVARIANT_MEASURE} suggest that $\{Z_N\}$ should converge to the solution of the SDE with instantaneous reflection at the boundary $z=0$, provided it is defined and unique. This property remains to be proved.
\medskip

The rest of this section is devoted to the proof of Theorem \ref{MAINRES} which, for the sake of clarity, is decomposed into three subsections.

\subsection{Proof of compactness}\label{S-COMPACT}
Recall that $\{P_N^\mathrm{x}\}_{\mathrm{x}\in S_N}$ denotes the Markov process generated by $L_N$. For every $\mathrm{x}\in S_N$, the probability measure $P_N^\mathrm{x}$ solves the martingale problem associated with $L_N$ and the initial condition $\mathrm{x}$. In particular the stochastic process $M_N$ on $(\D([0,T],S_N),{\cal F}_{T,N})$ with values in $\R$, defined by 
\[
t\mapsto M_N(t)=Z_N(t)-A_N(t)\, ,\quad \text{where}\quad A_N(t):=\int_0^t(L_Nz)(X_N(s))ds\, ,
\]
is a martingale relative to $P_N^{X_N(0)}$. Moreover, that $L_Nz$ is bounded in $S_N$ implies that $A_N$ is a process of finite variation. Consequently $M_N$ is bounded and then is locally square integrable. 

According to Corollary 2.3.3 in \cite{JM86}, in order to prove that the laws associated with $\{Z_N\}$ form a tight family, it suffices to show that 
\begin{itemize}
\item the sequences $\{A_N\}$ and $\{\langle M_N\rangle\}$ satisfy the Aldous condition and
\item the sequence of the laws of $\sup_{t\in [0,T]} |A_N(t)|$ (resp.\ $\sup_{t\in [0,T]} |\langle M_N(t)\rangle|$) is tight in $\R$. 
\end{itemize}
 Below we focus on proving the Aldous condition. The proof of the other condition is similar and left to the reader.

In order to prove the Aldous condition for $A_N$, we observe that the Markov inequality implies that for every $0\leq t< t'$, we have 
\[
P_N^\mathrm{x}(|A_N(t')-A_N(t)|\geq \eta)\leq \frac1{\eta}\E_\mathrm{x}\left(\int_{t}^{t'}|L_Nz(X_N(s))|ds\right),
\]
Explicit calculations using that ${\cal L}_{\rm fast}z=0$
imply that the integrand $|L_Nz(X_N(s))|$ is uniformly bounded in $N$, and so the probability $P_N^\mathrm{x}(|A_N(t')-A_N(t)|\geq \eta)$ can be made arbitrarily small, uniformly in $N$, by taking $t'-t$ sufficiently small.
 
Moreover, the increasing process $\langle M_N\rangle$ is given by 
\[
\langle M_N(t)\rangle = \int_0^t q_Nz(X_N(s))ds
\]
where $q_N$ is the quadratic operator defined for any function $f:S_N\to \R$ by 
\begin{equation}
q_Nf(\mathrm{x})=\sum_{i\in\Z_3}x_i(\frac{a}{N}+x_{i+1})N^2 \left(f(\mathrm{x}+\frac{\mathrm{u}_i}{N})-f(\mathrm{x})\right)^2.
\label{QUADOP}
\end{equation}
A similar argument as above applies to prove the Aldous condition for $\langle M_N\rangle$, using the mean value theorem for $z$.  

\subsection{Extending the time scale of the convergence to the deterministic approximation}
Let $\|\cdot\|_1$ denote the $\ell^1$-norm in $S$.
\begin{Lem}
Let $\{T_N\}_{N\in\N}$ be a sequence in $\R^+$ such that $T_N\leq C\frac{\log \log N}{N}$ for all $N\in \N$, for some $C>0$. Then we have 
\[
\lim_{N\to\infty}\sup_{t\in [0,T_N]}\sup_{X_N(0)\in S_N}\E_{P_N^{X_N(0)}} \left\| X_N(t)-X^{\rm fast}_{X_N(0)}(Nt)\right\|_1=0.
\]
\label{LIMITSHORTTIMES}
\end{Lem}

\noindent
{\sl Proof.} Given $f\in C^2(S)$, consider the martingale $M_N^f$ relative to $P_N^{X_N(0)}$ and defined by 
\begin{equation}
t\mapsto M_N^f(t):=f(X_N(t))-\int_0^tL_Nf(X_N(s))ds
\label{MARTF}
 \end{equation}
Using also the relation $f(X^{\rm fast}_{X_N(0)}(Nt))=f(X_{N}(0))+\int_0^tN{\cal L}_\text{fast}f(X^{\rm fast}_{X_N(0)}(s))ds$, we then get the estimate
\begin{align*}
\E_{P_N^{X_N(0)}} \left|f(X_N(t))-f(X^{\rm fast}_{X_N(0)}(Nt))\right|\leq & \E_{P_N^{X_N(0)}} \left|M_N^f(t)-f(X_{N}(0))\right|\\
&+\int_0^t\E_{P_N^{X_N(0)}} \left|L_Nf(X_N(s))-N{\cal L}_\text{fast}f(X_N(s))\right|ds\\
&+N\int_0^t\E_{P_N^{X_N(0)}} \left|{\cal L}_\text{fast}f(X_N(s))-{\cal L}_\text{fast}f(X^{\rm fast}_{X_N(0)}(s))\right|ds.
\end{align*}
We estimate each term in the RHS separately. By the H\"older inequality, we have 
\[
 \E_{P_N^{X_N(0)}} \left|M_N^f(t)-f(X_{N}(0))\right|\leq \left(\E_{P_N^{X_N(0)}} \left| M_N^f(t)-f(X_{N}(0))\right|^2\right)^{\tfrac12}= \left(\E_{P_N^{X_N(0)}} \langle M_N^f(t)\rangle\right)^{\tfrac12},
 \]
 where $\langle M_N^{f}(t)\rangle = \int_0^t q_Nf(X_N(s))ds$ and the quadratic operator $q_N$ is defined in \eqref{QUADOP}. This definition implies that for every $f\in C^1(S)$, there exists $K_f>0$ such that 
 \[
 \sup_{\mathrm{x}\in S,N\in\N}|q_Nf(\mathrm{x})|\leq K_f
 \]
 from where we get the upper bound $\E_{P_N^{X_N(0)}} \left|M_N^f(t)-f(X_{N}(0))\right|\leq K_f\sqrt{t}$.
 
Similarly, given $f\in C^2(S)$, the Taylor theorem applied to the first order expansion of $f(\mathrm{x}+\frac{\mathrm{u}_i}{N})$ at $\mathrm{x}$ implies the inequality
\[
\int_0^t\E_{P_N^{X_N(0)}} \left|L_Nf(X_N(s))-N{\cal L}_\text{fast}f(X_N(Ns))\right|ds\leq K_f t,
\]
provided that $K_f$ is chosen sufficiently large. By choosing $K_f$ even larger if necessary so that the following inequality holds 
\[
|{\cal L}_\text{fast}f(\mathrm{x})-{\cal L}_\text{fast}f(\mathrm{y})|\leq K_f\|\mathrm{x}-\mathrm{y}\|_1,\quad \forall \mathrm{x},\mathrm{y}\in S,
\]
and by collecting all the estimates, we finally obtain
\[
\E_{P_N^{X_N(0)}} \left|f(X_N(t))-f(X^{\rm fast}_{X_N(0)}(Nt))\right|\leq K_f\left(\sqrt{t}+t+N\int_0^tE_{P_N^{X_N(0)}} \left\|X_N(s)-X^{\rm fast}_{X_N(0)}(Ns)\right\|_1ds\right).
\]
Applying this inequality to the functions $f(\mathrm{x})=x_i$ ($i\in \Z_3$), and using the Gronwall inequality, we finally get
\[
\sup_{X_N(0)\in S_N}\E_{P_N^{X_N(0)}} \left\|X_N(t)-X^{\rm fast}_{X_N(0)}(Nt)\right\|_1\leq K\left(\sqrt{t}+t\right) e^{NKt},
\]
for some $K>0$. The Lemma then immediately follows from the fact that $T_N\leq C\frac{\log \log N}{N}$ implies that $\lim_{N\to\infty} (\sqrt{T_N}+T_N)e^{NKT_N}=0$. \hfill $\Box$

\subsection{Identification of the limits, end of the proof of Theorem \ref{MAINRES}}\label{S-INDENT}
The tightness of the sequence $\{Z_N\}$ (section \ref{S-COMPACT}) implies that, up to passing to a subsequence, this sequence converges in law to a process $Z$ with values in $\overline{I}$. To complete the proof of Theorem \ref{MAINRES}, it remains to prove that $Z$ must be a solution of the SDE \eqref{eq:eds_avg}. The core argument is to establish that $Z$ must solve the martingale problem associated with ${\cal L}_\text{avg}$ for a sufficiently large set of functions. To that goal, we first invoke a slight extension of the Skorokhod's representation theorem, see Theorem \ref{THMSKO} in Appendix \ref{A-SKO}, according to which there exist a common probability space $(\mathrm{E},\Omega,\mathrm{P})$ in which the random variables $z(X_N)$ pointwise converge to $Z$. Then, we consider the martingales $M_N^{f_I\circ z}$ defined by 
\[
t\mapsto M_N^{f_I\circ z}(t)=f_I(Z_N(t))-\int_0^tL_N(f_I\circ z)(X_N(s))ds.
\]
The key argument of the proof of Theorem \ref{MAINRES} is the following $L^1$-convergence of the martingales. 
\begin{Pro}
For every $f_I\in C^3(\overline{I})$ and $T>0$, we have
\[
\lim_{N\to \infty}\E_{\mathrm{P}}\left|M_N^{f_I\circ z}(t)-f_I(Z(t))-\int_0^t{\cal L}_\text{\rm avg}f_I(Z(s))ds\right|=0,\quad \forall t\in [0,T].
\]
\end{Pro}
\noindent
The end of the proof of the Theorem \ref{MAINRES} uses again standard arguments. The Proposition implies in particular that the process $t\mapsto f_I(Z(t))+\int_0^t{\cal L}_\text{\rm avg}f_I(Z(s))ds$ equipped with the probability $\mathrm{P}$ is a martingale for every $f_I\in C^3(\overline{I})$, see e.g.\ Lemma 3.6, Chap.\ 2 in \cite{D96}, and hence for $f_I(Z)=Z$ and $f_I(Z)=Z^2$. Moreover, the process is continuous (proved in the proof of the proposition below). Hence by (an adaptation to $\overline{I}$ of the) Proposition 4.6, Chap.\ 5 in \cite{KS91}, one defines a Brownian motion $W$ such that the process is a solution of the SDE \eqref{eq:eds_avg}.
\medskip

\noindent
{\sl Proof of the Proposition.} We are going to expand the difference $M_N^{f_I\circ z}(t)-f_I(Z(t))+\int_0^t{\cal L}_\text{avg}f_I(Z(s))ds$ into a telescopic sum for which the elements can be controlled using characteristic features of the various processes involved in the approximation. 

To that goal, let $k_N=\lfloor\tfrac{tN}{\log \log N}\rfloor$ and $T_N=\tfrac{t}{k_N}$. Using the definition of $M_N^{f_I\circ z}$ above, we write
\begin{equation}
\E_{\mathrm{P}}\left|M_N^{f_I\circ z}(t)-f_I(Z(t))+\int_0^t{\cal L}_\text{avg}f_I(Z(s))ds\right|\leq \E_{\mathrm{P}}\left|f_I(Z_N(t))-f_I(Z(t))\right|+\sum_{\ell=1}^5\E_{\mathrm{P}}\left|Q_{N,\ell}^{f_I}(t)-Q_{N,\ell+1}^{f_I}(t)\right|,
\label{INEQMART}
\end{equation}
where the boundary terms 
\[
Q_{N,1}^{f_I}(t)=\displaystyle\int_0^tL_N(f_I\circ z)(X_N(s))ds\quad \text{and}\quad Q_{N,6}^{f_I}(t)=\int_0^t{\cal L}_\text{avg}f_I(Z(s))ds
\] 
in the telescopic sum are already known (NB: While $Q_{N,6}^{f_I}(t)$ actually does not depend on $N$, using this notation simplifies the expression \eqref{INEQMART}) and the interior terms $Q_{N,\ell}^{f_I}(t)$ are defined by
\[
Q_{N,\ell}^{f_I}(t)=
\begin{cases}
{\displaystyle\int_0^t{\cal L}_\text{slow}(f_I\circ z)(X_N(s))ds}\quad\text{if}\quad \ell=2\\
{\displaystyle\sum_{k=0}^{k_N-1}\int_{kT_N}^{(k+1)T_N}{\cal L}_\text{slow}(f_I\circ z)(X^{\rm fast}_{X_N(kT_N)}(Ns))ds}\quad\text{if}\quad \ell=3\\
{\displaystyle T_N\sum_{k=0}^{k_N-1}{\cal L}_\text{avg}f_I(Z_N(kT_N))}\quad\text{if}\quad \ell=4\\
{\displaystyle T_N\sum_{k=0}^{k_N-1}{\cal L}_\text{avg}f_I(Z(kT_N))}\quad\text{if}\quad \ell=5\\
\end{cases}
 \]
 We now prove that each term in the RHS of \eqref{INEQMART} vanishes in the limit of large $N$.
 \medskip
 
\noindent
{\bf Proof of convergence of $\E_{\mathrm{P}}\left|f_I(Z_N(t))-f_I(Z(t))\right|$ and $\E_{\mathrm{P}}\left|Q_{N,4}^{f_I}(t)-Q_{N,5}^{f_I}(t)\right|$}. We first write
\[
\E_{\mathrm{P}}\left|f_I(Z_N(t))-f_I(Z(t))\right|\leq \E_{\mathrm{P}}\left\|f_I(Z_N)-f_I(Z)\right\|_\infty,\ t\in [0,T],
\]
where $\|\cdot\|_\infty$ denotes the uniform norm on $[0,T]$. Similarly, we have
\[
\E_{\mathrm{P}}\left|Q_{N,4}^{f_I}(t)-Q_{N,5}^{f_I}(t)\right|\leq t \E_{\mathrm{P}}\left\|{\cal L}_\text{avg}f_I(Z_N)-{\cal L}_\text{avg}f_I(Z)\right\|_\infty,\ t\in [0,T].
\]
The amplitude of the jumps of $Z_N$ is at most $\tfrac1{N}$ and the mapping $(Z(t))_{t\in\R^+}\mapsto \sup_{t\in [0,T]}|Z(t)-Z(t^-)|$ is continuous in the Skorokhod topology of $\D([0,T],\overline{I})$; hence the limit $Z$ must be continuous a.s.\ and the convergence $Z_N\to Z$ a.s.\ must occur in the sense of $\|\cdot\|_\infty$. The desired convergences then follow from the dominated convergence theorem using that $f_I$ and ${\cal L}_\text{avg}f_I$ are uniformly continuous over $\overline{I}$.
 \medskip
 
 \noindent
 {\bf Proof of convergence of $\E_{\mathrm{P}}\left|Q_{N,1}^{f_I}(t)-Q_{N,2}^{f_I}(t)\right|$}. Given $f\in C^3(S)$, the Taylor theorem applied to the second order expansion of $f(\mathrm{x}+\frac{\mathrm{u}_i}{N})$ at $\mathrm{x}$ implies the existence of $K_f>0$ such that we have 
 \[
 \left|L_Nf(\mathrm{x})-N{\cal L}_\text{fast}f(\mathrm{x})-{\cal L}_\text{slow}f(\mathrm{x})\right|\leq \frac{K_f}{N}.
 \]
Applying this inequality to $f_I\circ z$ with $f_I\in C^3(\overline{I})$, and using the property ${\cal L}_\text{fast}(f_I\circ z)=0$ (which follows from the fact that the function $\mathrm{z}$ is invariant under the flow generated by $F$), we immediately obtain the desired convergence
\[
\lim_{N\to\infty}\E_{\mathrm{P}}\left|Q_{N,1}^{f_I}(t)-Q_{N,2}^{f_I}(t)\right| =0,\ \forall t\in [0,T].
\]
\medskip
 
 \noindent
 {\bf Proof of convergence of $\E_{\mathrm{P}}\left|Q_{N,2}^{f_I}(t)-Q_{N,3}^{f_I}(t)\right|$}. Given $f_I\in C^2(\overline{I})$, let $K_s$ be the Lipschitz constant of ${\cal L}_\text{slow}(f_I\circ z)$ with respect to the $\|\cdot\|_1$-norm. We have
 \begin{align*}
 \E_{\mathrm{P}}\left|Q_{N,2}^{f_I}(t)-Q_{N,3}^{f_I}(t)\right|&\leq K_s\sum_{k=0}^{k_N-1}\int_{kT_N}^{(k+1)T_N}\E_{\mathrm{P}}\|X_N(s)-X^{\rm fast}_{X_N(kT_N)}(Ns)\|_1ds\\
 &\leq K_s t \sup_{t\in [0,T_N]}\sup_{X_N(0)\in S_N}\E_{P_N^{X_N(0)}} \left\| X_N(t)-X^{\rm fast}_{X_N(0)}(Nt)\right\|_1,
 \end{align*}
and then Lemma \ref{LIMITSHORTTIMES} immediately imply 
\[
\lim_{N\to\infty}\E_{\mathrm{P}}\left|Q_{N,2}^{f_I}(t)-Q_{N,3}^{f_I}(t)\right| =0,\ \forall t\in [0,T].
\] 
\medskip
 
 \noindent
 {\bf Proof of convergence of $\E_{\mathrm{P}}\left|Q_{N,3}^{f_I}(t)-Q_{N,4}^{f_I}(t)\right|$}. The proof of convergence for this term follows from considerations about localisation in $S$ and related dynamical estimates. Writing 
 \begin{equation}
 \E_{\mathrm{P}}\left|Q_{N,3}^{f_I}(t)-Q_{N,4}^{f_I}(t)\right|\leq T_N\sum_{k=0}^{k_N-1}\E_{\mathrm{P}}\left|\frac1{T_N}\int_{kT_N}^{(k+1)T_N}{\cal L}_\text{slow}(f_I\circ z)(X^{\rm fast}_{X_N(kT_N)}(Ns))ds-{\cal L}_\text{avg}f_I(Z_N(kT_N))\right|
 \label{ESTQ34}
 \end{equation}
 and given $r\in I$ (to be specified later on), for each term in the sum of the RHS, we consider separately the cases $Z_N(kT_N) \in [0,r)$ and $Z_N(kT_N)\in [r,\tfrac1{27}]$. 
 
 In the second case, we use that given $t>0$ and $\mathrm{x}\in S$ such that $z(\mathrm{x})\geq r$, $\mathrm{x}\neq (\tfrac13,\tfrac13,\tfrac13)$, the definition of the averaged generator ${\cal L}_\text{avg}$ and the periodicity $X^{\rm fast}_{\mathrm{x}}(s+T(z(\mathrm{x})))= X^{\rm fast}_{\mathrm{x}}(s)$ imply
 \begin{align*}
 &\left|\tfrac1{NT_N}\int_0^{NT_N}{\cal L}_\text{slow}(f_I\circ z)(X^{\rm fast}_{\mathrm{x}}(s))ds-{\cal L}_\text{avg}f_I(Z(\mathrm{x}))\right| \\
 &\leq  \tfrac1{NT_N}\left|\int_0^{NT_N}{\cal L}_\text{slow}(f_I\circ z)(X^{\rm fast}_{\mathrm{x}}(s))ds-\lfloor\frac{NT_N}{T(z(\mathrm{x}))}\rfloor\int_0^{T(z(\mathrm{x}))}{\cal L}_\text{slow}(f_I\circ z)(X^{\rm fast}_{\mathrm{x}}(s))ds\right|\\
 &+\left|\left(\frac{T(z(\mathrm{x}))}{NT_N}\lfloor\frac{NT_N}{T(z(\mathrm{x}))}\rfloor-1\right){\cal L}_\text{avg}f_I(Z(\mathrm{x}))\right|\\
 &\leq \tfrac1{NT_N}\left|\int_{T(z(\mathrm{x}))\lfloor\frac{NT_N}{T(z(\mathrm{x}))}\rfloor}^{NT_N}{\cal L}_\text{slow}(f_I\circ z)(X^{\rm fast}_{\mathrm{x}}(s))ds\right|+\frac{T(z(\mathrm{x}))}{NT_N}\left|{\cal L}_\text{avg}f_I(Z(\mathrm{x}))\right|.
 \end{align*}
Accordingly, and using also that the period $T(z(\mathrm{x}))$ is bounded over those points $\mathrm{x}\in S$ for which $z(\mathrm{x})\geq r$ (see Lemma \ref{PROPERIOD} {\em (i)}), we conclude that the RHS vanishes in the limit $N\to\infty$, since $\lim_{N\to\infty}NT_N=+\infty$. For $\mathrm{x}=(\tfrac13,\tfrac13,\tfrac13)$, the result immediately follows from the equality \eqref{MATCH}.
 
The first case $Z_N(kT_N) \in [0,r)$ corresponds to the neighborhood of the simplex boundaries, where the control of averaging is more elusive. Recall from the comments after Lemma \ref{AVGOP} that ${\cal L}_\text{avg}f_I(0^+)=0$. Hence by taking $r$ sufficiently small, for each (putative) term ${\cal L}_\text{avg}f_I(Z_N(kT_N))$ in \eqref{ESTQ34}, we can make its contribution arbitrarily small (and the same comment applies to the maximal total contribution $T_N\sum_{k=0}^{k_N-1}\left|{\cal L}_\text{avg}f_I(Z_N(kT_N))\right|$).
 
In order to address the remaining integral term in \eqref{ESTQ34}, we observe that when $Z_N(t) \in [0,r)$, we have $z(X^{\rm fast}_{X_N(t)}(Ns))\in [0,r)$ for all $s\in\R^+$. Therefore, at any $s$, the point $X^{\rm fast}_{X_N(t)}(Ns)$ may be close to one of the vertices $\mathrm{v}_i$ of $S$.  
 An explicit calculation shows that 
 \[
 \lim_{\mathrm{x}\to \mathrm{v}_i}{\cal L}_\text{slow}(f_I\circ z)(\mathrm{x})=0,\ i\in\Z_3.
 \]
 As before, to choose some neighbourhoods $V_{i,r}$ of the $\mathrm{v}_i$ with sufficiently small radius $r$ implies that the contribution of the integral terms in \eqref{ESTQ34}, for those $X^{\rm fast}_{X_N(kT_N)}(Ns)\in \bigcup_{i\in\Z_3}V_{i,r}$, can be made arbitrarily small, uniformly in $N$.
 
W.l.o.g.\ we may assume that $\bigcup_{i\in\Z_3}V_{i,r}$ is invariant under the cyclic permutation of coordinates $(x_i)\mapsto (x_{i+1})$. Then, if a trajectory $X^{\rm fast}_{X_N(t)}(Ns)$ leaves a set $V_{i,r}$, then it must travel to $V_{i+1,r}$. In the intermediate region between $V_{i,r}$ and $V_{i+1,r}$, the norm $\|F(\mathrm{x})\|$ of the vector field is bounded below (Indeed, one easily checks that this is the case when restricted to the segment $\mathrm{v}_i\mathrm{v}_{i+1}$, part of an edge in the boundary $\partial S$. Then apply a continuity argument). Therefore, the transit time of $X^{\rm fast}_{X_N(t)}(Ns)$ must be bounded from below by, say $\tfrac{t_r}{N}$ for some $t_r>0$. It follows that, in the interval $[kT_N,(k+1)T_N]$ the total time the trajectory spends outside $\bigcup_{i\in\Z_3}V_{i,r}$, cannot exceed $3T_N\tfrac{t_r}{N}$. The corresponding total contribution of the integral terms in \eqref{ESTQ34}, for those $X^{\rm fast}_{X_N(kT_N)}(Ns)$ in the transit regions between the $V_{i,r}$, then cannot exceed $3\tfrac{t_rk_N}{N}$, which vanishes when $NT_N\to\infty$. 
This completes the proof that 
\[
\lim_{N\to\infty}\E_{\mathrm{P}}\left|Q_{N,3}^{f_I}(t)-Q_{N,4}^{f_I}(t)\right| =0,\ \forall t\in [0,T].
\]
\medskip
 
 \noindent
 {\bf Proof of convergence of $\E_{\mathrm{P}}\left|Q_{N,5}^{f_I}(t)-Q_{N,6}^{f_I}(t)\right|$}. The quantity ${\displaystyle T_N\sum_{k=0}^{k_N-1}{\cal L}_\text{avg}f_I(Z(kT_N))}$ can be regarded as a Riemann sum for the integral $\int_0^t{\cal L}_\text{avg}f_I(Z(s))ds$; hence the desired convergence follows since ${\cal L}_\text{avg}f_I\circ Z$ in \textit{a.s} continuous over $[0,T]$.
 
\noindent
The proof of the Proposition is complete. \hfill $\Box$ 
 
 \medskip
 
 \noindent
 {\large \bf Acknowledgements} 
 
 \noindent
This work has been supported by the ANR-19-CE40-0023 (PERISTOCH). We are grateful to Nils Berglund, Nicolas Fournier and Luc Hillairet for fruitful discussions and relevant comments.

\appendix

\section{Proof of Proposition \ref{MU_INV_N_INFINI}}\label{A-MU_INV_N_INFINI}
One ingredient is the following convergence of the uniform atomic measure $\nu_N=\tfrac1{|S_N|}\sum_{\mathrm{x}\in S_N}\delta_{\mathrm{x}}$ on $S$ to the Lebesgue measure.
\begin{Claim}
We have $\lim_{N\to\infty}\nu_N=\text{Leb}_{S}$.
\label{CONVMULEB}
\end{Claim}

\noindent
{\sl Sketch of proof of the Claim:} This property is a consequence of the fact that every continuous function on $S$ is Riemann integrable on this set. Consider the restriction to $S$ of the Vorono\"{\i} triangulation of the points in $S_N$. Since the distribution of these points is uniform and isotropic, the cells associated with points in the interior of $S$ are all copies of the same polyhedron $P_N$, whose diameter vanishes as $N\to\infty$. Therefore, given $f\in C^0(S)$, the sum 
\[
\text{Vol}(P_N)\sum_{\mathrm{x}\in \text{Int}(S_N)}f(\mathrm{x})+\sum_{\mathrm{x}\in S_N\setminus \text{Int}(S_N)}\text{Vol}(P_N(\mathrm{x})\cap S)f(\mathrm{x})
\]
(where $P_N(\mathrm{x})$ is the cell at $\mathrm{x}\in  S_N\setminus \text{Int}(S_N)$) can be viewed as a Darboux sum for the integral $\int_Sf(\mathrm{x})d\mathrm{x}$. The claimed weak convergence then readily follows from the following estimates
\[
\frac1{|S_N|}<\frac{\text{Vol}(P_N)}{\text{Vol}(S)}<\frac1{|\text{Int}(S_N)|}\quad\text{and}\quad \lim_{N\to\infty}\frac{|S_N\setminus \text{Int}(S_N)|}{|S_N|}=0.
\]
\hfill $\Box$

Independently, the Wendel limit (see e.g.\ \cite{AS70}, page 257, 6.1.46) implies
\[
\lim_{N\to\infty}\frac1{N^{a-1}}\frac{\Gamma(Nx+a)}{\Gamma(Nx+1)}=x^{a-1},\ \forall x\geq 0\ \text{if}\ a\geq 1\ (\forall x>0\ \text{if}\ a< 1)
\]
which suggests to consider the convergence of densities 
\[
\lim_{N\to\infty}\rho_{N,a}(\mathrm{x})=\rho_a(\mathrm{x})\quad\text{where}\quad \rho_{N,a}(\mathrm{x})=\frac{1}{N^{3(a-1)}}\prod_{i\in\Z_3}\frac{\Gamma(Nx_i+a)}{\Gamma(Nx_i+1)}.
\]
An analysis of the sign of the derivative of the function $y\mapsto \ln\left(\tfrac{\Gamma (yx+a)}{y^{a-1}\Gamma(yx+1)}\right)$ for $x\geq 0$ yields the following conclusion
\begin{itemize}
\item[] If $a< 1$, the sequence $\left\{\tfrac1{N^{a-1}}\frac{\Gamma(Nx+a)}{\Gamma(Nx+1)}\right\}_{N\in\N}$ is increasing for all $x>0$.
\item[] If $a\geq 1$, the sequence $\left\{\tfrac1{N^{a-1}}\frac{\Gamma(Nx+a)}{\Gamma(Nx+1)}\right\}_{N\in\N}$ is non-increasing for all $x\geq 0$.
\end{itemize}
Dini's Theorem then implies that the convergence to the Wendel limit is uniform in every compact set included in $(0,1]$ when $a<1$ (resp.\ in $[0,1]$ when $a\geq 1$).

In order to prove that $\lim_{N\to\infty}\mu_{N,a}=\mu_a$, we separate the case $a\geq 1$ and $a\in (0,1)$. 
 
If $a\geq 1$, the proof is immediate. Indeed, given $f\in C^0(S)$, the uniform convergence above and Claim \ref{CONVMULEB} respectively imply
\[
\lim_{N\to\infty}\int_S f(\mathrm{x})\left(\rho_{N,a}(\mathrm{x})-\rho_a(\mathrm{x})\right)d\nu_N(\mathrm{x})=0
\ \text{and}\ \lim_{N\to\infty} \int_S f(\mathrm{x})\rho_a(\mathrm{x}) d\nu_N(\mathrm{x})=\int_S f(\mathrm{x})\rho_a(\mathrm{x}) d\mathrm{x}
\]
which immediately yields
\[
\lim_{N\to\infty}\int_S f(\mathrm{x})\rho_{N,a}(\mathrm{x})d\nu_N(\mathrm{x})=\int_S f(\mathrm{x})\rho_a(\mathrm{x}) d\mathrm{x},
\]
and then the desired convergence, considering that for $f(\mathrm{x})=1$, this previous relation gives $\lim_{N\to\infty}\left(C_{N,a}|S_N|N^{3(a-1)}\right)^{-1}=C_a^{-1}$.

For $a\in (0,1)$, the proof is more involved because $\rho_a$ is not defined on the boundary of $S$ while this set has positive measure for $\mu_{N,a}$ (hence the limitations on the domains where uniform convergence holds). Let again $f\in C^0(S)$ and given $\epsilon>0$ arbitrary, let $\delta>0$ be sufficiently small so that 
\[
\int_{S\setminus S(\delta)}\left|f(\mathrm{x})\right|\rho_a(\mathrm{x}) d\mathrm{x}<\tfrac{\epsilon}5\quad\text{where}\quad S(\delta)=\left\{\mathrm{x}\in S\ :\ x_i\geq \delta,\ \forall i\right\}.
\]
Besides, a similar decomposition as for $a\geq 1$ and uniform convergence on $S(\delta)$ imply that we have
\[
\left|\int_{S(\delta)} f(\mathrm{x})\rho_{N,a}(\mathrm{x})d\nu_N(\mathrm{x})-\int_{S(\delta)} f(\mathrm{x})\rho_a(\mathrm{x}) d\mathrm{x}\right|<\tfrac{\epsilon}5
\]
provided that $N$ is sufficiently large. It remains to consider the term $\int_{S\setminus S(\delta)}\left|f(\mathrm{x})\right|\rho_{N,a}(\mathrm{x}) d\nu_N(\mathrm{x})$, which we separate into two integrals, using the decomposition 
\[
S\setminus S(\delta)=(S\setminus \text{Int}(S))\cup (\text{Int}(S)\setminus S(\delta)).
\]
On one hand, the inequality $\rho_{N,a}\leq \rho_a$ and Claim \ref{CONVMULEB} imply
\[
\int_{\text{Int}(S)\setminus S(\delta)}\left|f(\mathrm{x})\right|\rho_{N,a}(\mathrm{x}) d\nu_N(\mathrm{x})\leq \int_{\text{Int}(S)\setminus S(\delta)}\left|f(\mathrm{x})\right|\rho_a(\mathrm{x}) d\nu_N(\mathrm{x})\leq \int_{S\setminus S(\delta)}\left|f(\mathrm{x})\right|\rho_a(\mathrm{x}) d\mathrm{x}+\tfrac{\epsilon}5\leq \tfrac{2\epsilon}5
\]
provided that $N$ is sufficiently large. On the other hand, control of the integral 
\[
\int_{S\setminus \text{Int}(S)}\left|f(\mathrm{x})\right|\rho_{N,a}(\mathrm{x}) d\nu_N(\mathrm{x})
\]
is provided by the following property.
\begin{Claim}
We have $\lim_{N\to\infty}\int_{S\setminus \text{Int}(S)}\left|f(\mathrm{x})\right|\rho_{N,a}(\mathrm{x}) d\nu_N(\mathrm{x})=0$.
\end{Claim}

\noindent
{\sl Proof of the Claim:} Since $f$ is bounded, it suffices to prove the result for $f(\mathrm{x})=1$. To that goal, the objects above are considered in arbitrary dimension (and the explicit dependence on $a$ is removed). Namely, let $D\in\N$, $\Z_D= \Z/D\Z$ and
\[
S_D=\left\{(x_i)_{i\in\Z_D}\in (\R^+)^{\Z_D}\ \text{such that}\ \sum_{i\in Z_D}x_i= 1\right\},\quad S_{N,D}=S_D\cap \tfrac1{N}\N^D,
\]
and 
\[
\rho_{D}(\mathrm{x})=C_D\left(\prod_{i\in\Z_D}x_i\right)^{a-1},\quad \rho_{N,D}(\mathrm{x})=\frac{1}{N^{D(a-1)}}\prod_{i\in\Z_D}\frac{\Gamma(Nx_i+a)}{\Gamma(Nx_i+1)}
\]
and $\nu_{N,D}=\tfrac1{|S_{N,D}|}\sum_{\mathrm{x}\in S_{N,D}}\delta_{\mathrm{x}}$. Now, consider the integral $I_{N,D}$ defined by $I_{N,D}=\int_{S_D\setminus \text{Int}(S_D)}\rho_{N,D}(\mathrm{x}) d\nu_{N,D}(\mathrm{x})$. The permutation symmetry implies that we have
\begin{align*}
I_{N,D}&\leq D\int_{S_D\cap\{\mathrm{x}: x_D=0\}}\rho_{N,D}(\mathrm{x}) d\nu_{N,D}(\mathrm{x})\\
&=D\frac{|S_{N,D-1}|}{|S_{N,D}|}\frac{\Gamma(a)}{N^{a-1}}\int_{S_{D-1}}\rho_{N,D-1}(\mathrm{x}) d\nu_{N,D-1}(\mathrm{x})\\
&=D\frac{D(D-1)N}{N+D-1}\frac{\Gamma(a)}{N^{a}}\left(I_{N,D-1}+\int_{\text{Int}(S_{D-1})}\rho_{N,D-1}(\mathrm{x}) d\nu_{N,D-1}(\mathrm{x})\right)\\
&\leq \frac{D(D-1)N}{N+D-1}\frac{\Gamma(a)}{N^{a}}\left(I_{N,D-1}+\frac2{C_{D-1}}\right),
\end{align*}
provided that $N$ is sufficiently large, where we used
%\footnote{The expression of $|S_{N,D}|$ can be obtained from the induction $|S_{N,D}|=\sum_{n=0}^N|S_{N,D-1}|$ and $|S_{N,1}|=1$.}
$|S_{N,D}|={N+D-1\choose D-1}$, the inequality $\rho_{N,D-1}\leq\rho_{D-1}$, Claim \ref{CONVMULEB} and $\int_{S_{D-1}}\rho_{D-1}(\mathrm{x})d\mathrm{x}=\tfrac1{C_{D-1}}$. The limit $\lim_{N\to\infty}I_{N,D}=0$ then easily follows by iterating this inequality over the decreasing values of $D$ and using also that $a>0$.
\hfill $\Box$

%\section{Sketch of proof of Proposition \ref{APPROXFAST}}\label{A-APPROXFAST}

\section{Analysis of the dynamics $\dot{\mathrm{x}}=F(\mathrm{x})$}\label{A-DYNAML0}
\subsection{Expressions of the solutions}\label{A-DYNAML0-1}
The equation $\dot{\mathrm{x}}=F(\mathrm{x})$ for $\mathrm{x}\in S$ gives the following system of two coupled ODEs 
\[
\left\{\begin{array}{l}
\dot{x_1}=x_1(1-x_1-2x_2)\\
\dot{x_2}=x_2(2x_1+x_2-1)
\end{array}\right.
\]
This system indicates that $x_1$ must increase when $x_2< x_3$ and must decrease when $x_3<x_2$. Similar considerations apply to $x_2$. Solving the equation $x_1 x (1-x-x_1)=z$ for $x$ (assuming $x_1\neq 0$) yields that we must have  
\[
x_2(x_1)=\tfrac{1-x_1-\sqrt{(1-x_1)^2-\frac{4z}{x_1}}}2\quad\text{and}\quad x_3(x_1)=\tfrac{1-x_1+\sqrt{(1-x_1)^2-\frac{4z}{x_1}}}2
\]
when $x_1$ increases and the RHS are exchanged when $x_1$ decreases. In particular, when $x_1$ increases, it must satisfy the following ODE 
\[
\dot{x_1}=x_1\sqrt{(1-x_1)^2-\frac{4z}{x_1}}.
\]
Moreover, by continuity, those values of $x\in (0,1)$ for which $x_2(x)=x_3(x)$ must be the roots of the polynomial $x(1-x)^2-4z=0$. It is direct to show that such roots exist and there are two of them, say $x^\text{min}(z)<x^\text{max}(z)$, iff $z\in I$. Alternatively $x^\text{min}(z)=\min_{t\in\R^+}x_1(t)$ and $x^\text{max}(z)=\max_{t\in\R^+}x_1(t)$. An explicit computation of the roots of the cubic polynomial yields the following expressions
\begin{eqnarray*}
x^\text{min}(z)=\tfrac23\left(1-\sin\left(\tfrac{\theta(27z)}{3}+\tfrac{\pi}6\right)\right)\\
x^\text{max}(z)=\tfrac23\left(1+\sin\left(\tfrac{\theta(27z)}{3}-\tfrac{\pi}6\right)\right)
\end{eqnarray*}
(and we have $x^\text{add}(z)=\frac23(1+\cos\tfrac{\theta(27z)}{3})$ for the third root) where the angle $\theta(z)$ is defined by the relations 
\[
\left\{\begin{array}{l}
\cos \theta(z)=2z-1\\
\sin \theta(z)=2\sqrt{z(1-z)}
\end{array}\right.
\]
In particular, the function $z\mapsto \theta(z)$ is decreasing over $(0,1)$ with range $(0,\pi)$ (and $\theta(1-z)=\pi-\theta(z)$), which implies the following limits 
\[
x^\text{min}(0^+)=0,\ x^\text{max}(0^+)=1\quad\text{and}\quad x^\text{min}(\tfrac1{27}^-)=x^\text{max}(\tfrac1{27}^-)=\tfrac13,
\]
as expected. Also, the function $z\mapsto \theta(z)$ is $C^\infty$ over $(0,1)$ (and its derivative is equal to $-\frac{1}{\sqrt{z(1-z)}}$). 

The dynamics $\dot{\mathrm{x}}=F(\mathrm{x})$ commutes with cyclic permutations of coordinates $x_i\mapsto x_{i+1}$. The permutation do not affect the product $x_1x_2x_3$; hence invariant loops are also invariant under these permutations. This implies the existence of $T>0$ such that $x_i(\cdot +T)=x_{i+1}(\cdot)$ for all $i$. Therefore, all trajectories must all be periodic. 

Moreover, the autonomous equation for $x_1$ in the previous subsection implies that  the half-period $\frac{T(z)}2$ can be defined as the time it takes for $x_1$ to transit from $x^\text{min}(z)$ to $x^\text{max}(z)$, namely we have
\[
T(z)=2\int_{x^\text{min}(z)}^{x^\text{max}(z)}\tfrac{dx}{\sqrt{x^2(1-x)^2-4zx}}.
\]

\subsection{Proof of Lemma \ref{PROPERIOD}}\label{A-PROPERIOD}
{\sl Proof that the function $z\mapsto T(z)$ is $C^\infty$ on $I$.} Given an initial point $(x_{1}^\text{in},x_{2}^\text{in})$, let $t\mapsto x_2(t,x_{1}^\text{in},x_{2}^\text{in})$ be the value of the second coordinate of the solution of the above system of coupled ODE's. Using that $x_2=\tfrac{1-x^\text{min}(z)}2$ (resp.\ $x_2=\tfrac{1-x^\text{max}(z)}2$) when $x_1$ reaches $x^\text{min}(z)$ (resp.\ $x^\text{max}(z)$), the period $T(z)$ can be (also) specified using the relation
\[
x_2(\tfrac{T(z)}2,x^\text{min}(z),\tfrac{1-x^\text{min}(z)}2)=\tfrac{1-x^\text{max}(z)}2.
\]
This relation can be regarded as an equation for $T(z)$ given $z$. In fact, the function $x_2(\frac{t}2,x_{1}^\text{in},x_{2}^\text{in})$ is $C^\infty$ jointly in $t$ and $(x_{1}^\text{in},x_{2}^\text{in})$, because the vector field itself is $C^\infty$. The functions $z\mapsto x^{\text{min}}(z)$ and $z\mapsto x^{\text{max}}(z)$ are $C^\infty$ on $I$ (inherited from the same property for $\theta(z)$). Also, for every $z\in I$, the derivative 
\[
\frac{d}{dt}x_2(\tfrac{T(z)}2,x^\text{min}(z),\tfrac{1-x^\text{min}(z)}2)=\dot x_2(x^\text{max}(z),\tfrac{1-x^\text{max}(z)}2)=\tfrac{1-x^\text{max}(z)}2\left(x^\text{max}(z)-\tfrac{1-x^\text{max}(z)}2\right)\neq 0,
\]
does not vanish. That $z\mapsto T(z)$ is $C^\infty$ then follows from the implicit function theorem.
\medskip

\noindent
{\sl Proof that $T(\frac1{27}^-)$ exists and is finite.} We observe that, together with the explicit expression of $T(z)$ above, the equality $x(1-x)^2-4z=(x-x^\text{min}(z))(x^\text{max}(z)-x)(x^\text{add}(z)-x)$ readily implies 
\[
\tfrac{2}{\sqrt{\overline{\chi}(z)}}\int_0^1\tfrac{dx}{\sqrt{x(1-x)}}\leq T(z)\leq \tfrac{2}{\sqrt{\underline{\chi}(z)}}\int_0^1\tfrac{dx}{\sqrt{x(1-x)}},
\]
where 
\[
\underline{\chi}(z)=\min_{x\in (x^\text{min}(z),x^\text{max}(z))}(x^\text{add}(z)-x)x\quad\text{and}\quad \overline{\chi}(z)=\max_{x\in (x^\text{min}(z),x^\text{max}(z))}(x^\text{add}(z)-x)x.
\]
We have the limit $x^\text{add}(\frac1{27}^-)=\frac43$ which implies the following ones $\underline{\chi}(\frac1{27}^-)=\overline{\chi}(\frac1{27}^-)=\frac13$. Together with $\int_0^1\frac{dx}{\sqrt{x(1-x)}}=\pi$, these behaviours yield the limit $T(\frac1{27}^-)=2\pi\sqrt{3}$.
%
%
%\noindent
%We write $z=\frac{1}{27}(1-\varepsilon)$. From their explicit expressions, we have the following expansions:
%\begin{equation}
%x^{\rm min}(z) =\frac13 -\frac{2}{3\sqrt{3}} \sqrt{\varepsilon} +O(\varepsilon);\quad x^{\rm max}(z) =\frac13 +\frac{2}{3\sqrt{3}} \sqrt{\varepsilon} +O(\varepsilon);\quad
%x^{\rm add}(z) =\frac43 +O(\varepsilon).
%\label{eq:dvpt_xmax127}
%\end{equation}
%We then factorize the polynomial expression: $x^2(1-x)^2-4xz = x(x-x^{\rm min}(z))(x^{\rm max}(z)-x)(x^{\rm add}(z)-x)$, and perform the change of variable
%$u=\frac{x-\frac13}{x^{\rm max}-\frac13}$.
%It is now a simple computation:
%\begin{eqnarray*}
%T(z) &=& 2   \int_{x^{\rm min}(z)}^{x_{\rm max}(z)} \frac{dx}{\sqrt{x(x-x^{\rm min}(z))(x^{\rm max}(z)-x)(x^{\rm add}(z)-x)}} \\
%&=& 2\frac{(x^{\rm max}(z)-\frac13)}{\sqrt{(\frac13-x^{\rm min}(z))(x^{\rm max}(z)-\frac13)}} \times \\
%&& \int_{\frac{x^{\rm min}(z)-\frac13}{x^{\rm max}(z)-\frac13}}^{1} \frac{du}{\left[\frac13+u(x^{\rm max}(z)-\frac13)\right]^{\frac12}\left[(x^{\rm add}(z)-\frac13-u(x^{\rm max}(z)-\frac13)\right]^{\frac12}\left[1-u\right]^{\frac12}\left[\frac{\frac13-x^{\rm min}(z)}{x^{\rm max}(z)-\frac13}+u\right]^{\frac12}}.
%\end{eqnarray*}
%Using the expansions \eqref{eq:dvpt_xmax127} and the dominated convergence theorem for $\varepsilon \to 0$, we conclude that $T(z)$ has a finite limit when $z$ tends to $\frac{1}{27}$.\hfill $\Box$
\medskip
\begin{figure}[ht]
\centerline{\includegraphics[width=7cm,trim={10cm 4cm 22cm 4cm},clip]{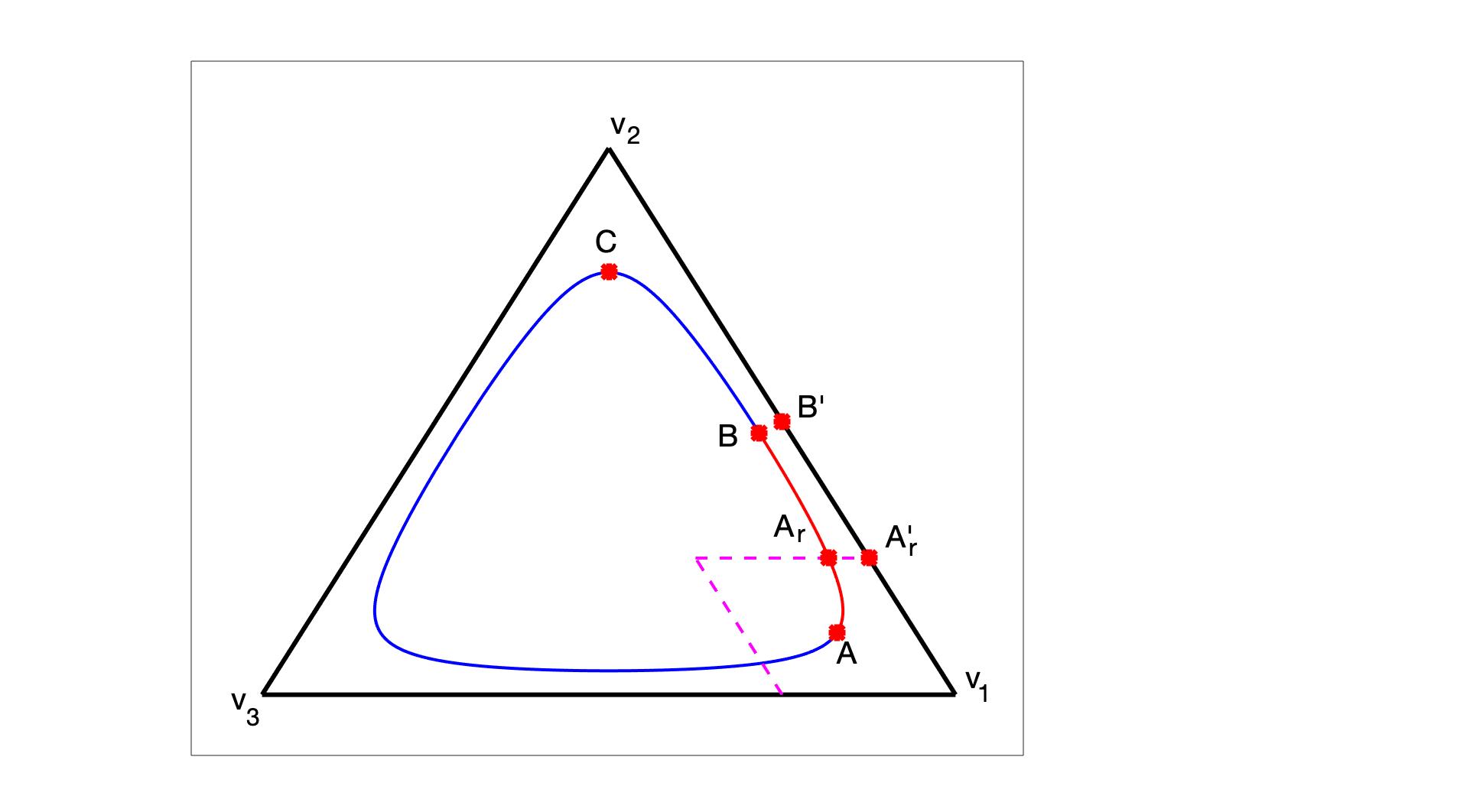}}
\caption{ An example for $z=0.01$ of periodic trajectory of the dynamics $\dot{\mathrm{x}}=F(\mathrm{x})$ (thin (blue) curve). The trajectory can be decomposed into 6 arcs, all of the same duration and which can be deduced one another using cyclic permutations of coordinates and/or time reversal; from $A$ (defined as the point where $x_1= x^{\rm max}(z)$) to $B$ (defined as $x_1=x_2=\tfrac{1-x^\text{min}(z)}2$), then to $C$ ($x_2=x^{\rm max}(z)$), etc. The (purple) dashed lines materialize the neighborhood $V_{r}$ of the vertex $\mathrm{v}_1$, such that ${\rm max}(x_2,x_3=1-x_1-x_2)\leq r$. The point $A_r$ is where the trajectory leaves $V_{r}$ (hence the coordinate $x_2$ of $A_r$ is $r$) and the segment $A'_rB'$ on the edge $\mathrm{v}_1\mathrm{v}_2$ is the segment that the arc $A_rB$ approaches as $z\to 0^+$. \label{fig:period}}
\end{figure}

\noindent
{\sl Proof that $T(z) \sim -3\ln z$ when $z\to 0^+$.} We decompose any periodic trajectory into 6 arcs, all of the same duration and which can be deduced one another using cyclic permutations of coordinates and/or time reversal, see Fig.\ \ref{fig:period}. Accordingly, let $\tau(z)=\frac{T(z)}6$ be the transit time in each arc. In particular, $\tau(z)$ is the transit time from $A$ (where $x_1= x^{\rm max}(z)$) to $B$ ($x_1=x_2=\tfrac{1-x^\text{min}(z)}2$). We aim at finding an equivalent for $\tau(z)$ as $z$ tends to $0^+$.

To that goal, let $r\in (0,\frac{1}{3})$ and consider the neighborhood  of the vertex $\mathrm{v}_1=(1,0)$, defined as $V_r=\{(x_1,x_2): {\rm max}(x_2,1-x_1-x_2)\leq r\}$. Let $z\in I$ be sufficiently small so that $A\in V_r$. Let then $A_r$ be the intersection point of the arc from $A$ to $B$ and the boundary of $V_r$. Let $\tau_{A\to A_r}$ (resp.\ $\tau_{A_r\to B}$) denote the transit time from $A$ to $A_r$ (resp.\ from $A_r$ to $B$). We clearly have 
\[
\tau(z)=\tau_{A\to A_r}+\tau_{A_r\to B}.
\]
We now estimate the two durations separately. For the latter, we have 
\[
\lim_{z\to 0^+}A_r=A'_r\quad\text{and}\quad \lim_{z\to 0^+} B= B',
\]
and the continuity of the vector field $F$ implies that $\lim_{z\to 0^+}\tau_{A_r\to B}=\tau_{A'_r\to B}$ which is finite for every $r>0$ (and diverges as $r\to 0^+$).

For the former, we observe that inside $V_r$, and hence along the arc from $A$ to $A_r$, the coordinate $x_1$ is not smaller than $1-2r$. Using also that $x_2\geq 0$, we obtain from the equation of the dynamics, the following inequalities for the time derivative of the coordinate $x_2$ in the arc between $A$ and $A_r$
\[
1-3r\leq \tfrac{\dot{x}_2}{x_2}\leq 1,
\]
from where direct integration yields the estimates
\[
(1-3r) \tau_{A\to A_r}\leq \ln x_2(A_r)-\ln\tfrac{1-x^{\text{\rm max}}(z)}2\leq \tau_{A\to A_r},
\]
and then
\[
1+\frac{\ln x_2(A_r)+\tau_{A_r\to B}}{-\ln\tfrac{1-x^{\text{\rm max}}(z)}2}\leq \frac{\tau(z)}{-\ln\tfrac{1-x^{\text{\rm max}}(z)}2}\leq \tfrac1{1-3r}+\frac{\frac{\ln x_2(A_r)}{1-3r}+\tau_{A_r\to B}}{-\ln\tfrac{1-x^{\text{\rm max}}(z)}2}.
\]
Now use that $x^{\text{\rm max}}(z)\sim 1-2\sqrt{z}$ as $z\to 0^+$, $\lim_{z\to 0^+} x_2(A_r)=x_2(A'_r)=r\sin\frac{\pi}3$ and the limit $\lim_{z\to 0^+}\tau_{A_r\to B}$ above to obtain
\[
1\leq \liminf_{z\to 0^+}\frac{\tau(z)}{-\frac12\ln z}\leq \limsup_{z\to 0^+}\frac{\tau(z)}{-\frac12\ln z}\leq\tfrac1{1-3r}.
\]
The desired conclusion finally follows by taking the limit $r\to 0$.

\subsection{Proof of Lemma \ref{pro:dyn_rapide}}\label{DYNRAP-A}
When evaluated for the trajectory such that $x_1(0)=x^\text{min}(z)$, using the change of variable $t\mapsto x_1(t)$, the quantity $A(z)$ can be written as follows
\begin{align*}
A(z)&=\int_{x^\text{min}(z)}^{x^\text{max}(z)}\tfrac{1-x_1-\sqrt{(1-x_1)^2-\frac{4z}{x_1}}}2dx_1+\int_{x^\text{max}(z)}^{x^\text{min}(z)}\tfrac{1-x_1+\sqrt{(1-x_1)^2-\frac{4z}{x_1}}}2dx_1\\
&=-\int_{x^\text{min}(z)}^{x^\text{max}(z)}\sqrt{(1-x_1)^2-\frac{4z}{x_1}}dx_1.
\end{align*}
The smoothness of the integration bounds and of the integrand imply, via the Leibniz integral rule, that the function $z\mapsto A(z)$ is $C^\infty$ and $I'(z)=T(z)$ on $I$. 

From the expression above, $A(z)$ can be interpreted as the area enclosed in the loop $z(\mathrm{x})=z$. Accordingly, we clearly have $A(0^+)=\frac12$ and $A(\frac1{27}^-)=0$. From the relation $A'(z)=T(z)$ on $I$ and $T(\frac1{27}^-)=2\pi\sqrt{3}$, we conclude that $A(z) \sim -2\pi\sqrt{3}(\frac1{27} - z)$ when $z \to \frac{1}{27}^-$.

\section{A slight extension of the Skorokhod's representation theorem}\label{A-SKO}
\begin{Thm}
Let $(E,{\cal B},P)$ and $(E_N,{\cal B}_N,P_N)$ for $N\in\N$ be probability spaces, where $E$ and $E_N$ are metric spaces, ${\cal B}$ and ${\cal B}_N$ their respective Borel $\sigma$-algebras and where the support of $P$ is separable. Assume the existence of measurable functions $f_N:E_N\to E$ such that $P_N\circ f_N^{-1}\Rightarrow P$ as $N\to \infty$ (weak convergence). Then there exist random variables $Y$ and $Y_N$ defined on a common probability space $(\mathrm{E},\Omega,\mathrm{P})$ such that the laws of $Y$ and $Y_N$ are respectively $P$ and $P_N$ and we have 
\[
\lim_{N\to \infty}Y_N(\omega)=Y(\omega)\ \text{for every}\ \omega\in \mathrm{E}.
\] 
\label{THMSKO}
\end{Thm}

\noindent
The proof is a direct adaptation {\sl mutatis mutandis} of the proof of Theorem 6.7, Chap.\ 1 in \cite{B99}.
\end{document}